\documentstyle[11pt]{article}
\newcommand{\nl}{ {\hfill \break} }
\newcommand{\np}{ {\newpage } }

\newcommand{\sm}{ {\setminus } }

\newcommand{\Iff}{ {\Leftrightarrow } }
\newcommand{\lto}{ {\longrightarrow } }
\newcommand{\imp}{ {\Rightarrow } }

\newcommand{\lam}{ {\lambda } }
\newcommand{\eps}{ \mbox{$\varepsilon$} }
\newcommand{\cl   }{ \mbox{${\rm cl   }$} }

\newcommand{\trace}{ \mbox{${\rm tr   }$} }

\newcommand{\N    }{ \mbox{${\rm I\!N    }$} }
\newcommand{\R    }{ \mbox{${\rm I\!R    }$} }
\newcommand{\GL   }{ \mbox{${\rm GL   }$} }

\newcommand{\NJNF }{ \mbox{${\rm NJNF }$} }
\newcommand{\ONJNF}{ \mbox{${\rm ONJNF}$} }
\newcommand{\I   }{ \mbox{${\rm I   }$} }
\newcommand{\II  }{ \mbox{${\rm II  }$} }
\newcommand{\III }{ \mbox{${\rm III }$} }
\newcommand{\IV  }{ \mbox{${\rm IV  }$} }
\newcommand{\Vv  }{ \mbox{${\rm V^{(4)} }$} }
\newcommand{\V   }{ \mbox{${\rm V   }$} }
\newcommand{\VI  }{ \mbox{${\rm VI}  $} }
\newcommand{\VIo }{ \mbox{${\rm VI}_0$} }
\newcommand{\VIh }{ \mbox{${\rm VI}_h$} }
\newcommand{\VII }{ \mbox{${\rm VII} $} }
\newcommand{\VIIo}{ \mbox{${\rm VII}_0$} }
\newcommand{\VIIh}{ \mbox{${\rm VII}_h$} }
\newcommand{\VIII}{ \mbox{${\rm VIII}$} }
\newcommand{\IX  }{ \mbox{${\rm IX  }$} }
\newcommand{\ii  }{ \mbox{${\rm ii  }$} }
\newcommand{\iv  }{ \mbox{${\rm iv  }$} }
\newcommand{\ve  }{ \mbox{${\rm v   }$} }
\newcommand{\IIRL  }{ \mbox{${\rm II  }^{R/L} $} }
\newcommand{\IVRL  }{ \mbox{${\rm IV  }^{R/L} $} }
\newcommand{\VIIoRL}{ \mbox{${\rm VII}^{R/L}_0 $} }
\newcommand{\VIIhRL}{ \mbox{${\rm VII}^{R/L}_h $} }
\newcommand{\VIIIRL}{ \mbox{${\rm VIII}^{R/L} $} }
\newcommand{\IXRL  }{ \mbox{${\rm IX  }^{R/L} $} }
\newcommand{\PI    }{ \mbox{${\wp_1   }$} }
\newcommand{\PII   }{ \mbox{${\wp_2   }$} }
\newcommand{\PIII  }{ \mbox{${\wp_3   }$} }
\newcommand{\PIV   }{ \mbox{${\wp_4   }$} }
\newcommand{\PV    }{ \mbox{${\wp_5   }$} }
\newcommand{\PVI   }{ \mbox{${\wp_6   }$} }
\newcommand{\PVIac }{ \mbox{${\wp^{\alpha/\gamma}_6 }$} }
\newcommand{\PVIa  }{ \mbox{${\wp^\alpha_6  }$} }
\newcommand{\PVIb  }{ \mbox{${\wp^\beta_6   }$} }
\newcommand{\PVIc  }{ \mbox{${\wp^\gamma_6  }$} }
\newcommand{\PVII  }{ \mbox{${\wp_7  }$} }
\newcommand{\PVIII }{ \mbox{${\wp_8  }$} }

\newcommand{\An    }{[1]}%
\newcommand{\Bi    }{[2]}%
\newcommand{\bBi   }{[3]}%
\newcommand{\Bo   }{[4]}
\newcommand{\Cha   }{[5]}%
\newcommand{\Ch   }{[6]}%
\newcommand{\Co   }{[7]}%
\newcommand{\Est  }{[8]}
\newcommand{\Gri  }{[9]}%
\newcommand{\Gru  }{[10]}%
\newcommand{\Hud  }{[11]}%
\newcommand{\In   }{[12]}%
\newcommand{\Ki   }{[13]}%
\newcommand{\Kr   }{[14]}%
\newcommand{\Lan  }{[15]}%
\newcommand{\Lee  }{[16]}%
\newcommand{\Lie  }{[17]}%
\newcommand{\bLie }{[18]}%
\newcommand{\Mac  }{[19]}%
\newcommand{\Mag  }{[20]}%
\newcommand{\Mo   }{[21]}%
\newcommand{\aMu  }{[22]}%
\newcommand{\bMu  }{[23]}%
\newcommand{\cMu  }{[24]}%
\newcommand{\Ne   }{[25]}%
\newcommand{\PaSWZ}{[26]}%
\newcommand{\PaW  }{[27]}%
\newcommand{\Pe   }{[28]}%
\newcommand{\aRa  }{[29]}
\newcommand{\bRa  }{[30]}
\newcommand{\Ri   }{[31]}
\newcommand{\Ro   }{[32]}%
\newcommand{\Sal  }{[33]}%
\newcommand{\San  }{[34]}%
\newcommand{\aSch }{[35]}%
\newcommand{\bSch }{[36]}%
\newcommand{\Se   }{[37]}%
\newcommand{\Smrz }{[38]}
\newcommand{\aTur }{[39]}
\newcommand{\bTur }{[40]}
\newcommand{\cTur }{[41]}
\newcommand{\Um   }{[42]}
\newcommand{\Vra  }{[43]}

\begin{document}

\begin{tabbing}
\` {\small UniP-Math-95/03}
  \\
\` {\small\it June 24, 1995}
\end{tabbing}

\vspace{1.1truecm}
\begin{center}
{{\bf
TOPOLOGICAL CLASSIFYING SPACES OF LIE ALGEBRAS
\vspace{0.2cm}\\
AND THE NATURAL COMPLETION OF CONTRACTIONS
\vspace{0.5cm}
}}
\end{center}
\begin{center}
{\bf {Martin Rainer}}
\dag \vspace{1.0truecm}\\
{Projektgruppe Kosmologie, Institut f\"ur Mathematik}\\
{Universit\"at Potsdam, Am Neuen Palais 10}\\
{PF 601553, D-14415 Potsdam, Germany}
\vspace{1.0truecm}
\end{center}
\begin{abstract}
\vspace{0.5truecm}
The space $K^n$ of all $n$-dimensional { Lie} algebras
has a natural non-Hausdorff topology ${\kappa}^n$, which
has characteristic limits, called transitions $A\to B$,
between distinct { Lie} algebras $A$ and $B$.
The entity of these transitions are the natural transitive completion
of the well known { In\"on\"u-Wigner} contractions and their
partial generalizations by { Saletan}.

Algebras containing a common ideal of codimension $1$
can be characterized by homothetically normalized { Jordan}
normal forms of one generator of their adjoint representation.
For such algebras, transitions $A\to B$ can be described
by limit transitions between corresponding normal forms.

The topology  $\kappa^n$ is presented in detail for $n\leq 4$.
Regarding the orientation of the algebras as vector spaces has a
non-trivial effect for the corresponding topological space $K^n_{or}$:
There exist both, selfdual points and pairs of dual points
w.r.t. orientation reflection.
\end{abstract}
\vspace{0.7truecm}
{\footnotesize
Financially supported by DFG grants Bl 365/1-1 and Schm 911/6-2.
}
\vspace{.17truecm}
\nl
{\small
\dag E-mail: mrainer@aip.de
}
\np
\section{\bf Introduction}
\setcounter{equation}{0}
The main goal of this paper is to study the topological space
of real {Lie} algebras of a given dimension $n\leq 4$.

Extensive studies have been dedicated to generalizations
of the classical {Lie} algebra structure. As an example think
of the famous q-deformations or Santilli's Lie isotopic
liftings \San. However,
few work has been dedicated to pursue the theory of
deformations and contractions of {Lie} algebras (or groups)
within their category.

{Smrz} \Smrz\ has considered the deformation of {Lie}
algebras outside a fixed subgroup. This kind of deformation
is in some sense complementary to a {In\"on\"u-Wigner} contraction
\In, which consists in a parametric linear and isotropic contraction
outside a given subgroup of a {Lie} algebra.

A particularly interesting problem is to find all possible
contractions and, more generally, all possible limit transitions
between real or complex {Lie} algebras of fixed
dimension $n$, and to uncover
the natural topological structure
of the space of all such {Lie} algebras.
It is clear that this requires, as a precondition,
to find all isomorphism classes of {Lie} algebras in the given
dimension.
Unfortunately, with increasing dimension $n$ the classification of
real and complex {Lie} algebras
becomes rapidly more complicated.

For this goal, the {Levi} decomposition into a semidirect
sum of a radical and a semisimple subalgebra proves to be useful.
This way {Turkowski} has classified
real {Lie} algebras which admit a nontrivial {Levi}
decomposition,  up to $n=8$ in \aTur\
and recently for $n=9$ in \cTur.

In any dimension $n$, the classification of all
nilpotent {Lie} algebras is an essential step required for a complete
classification.
For $n=7$, a complete list of all nilpotent, real and complex,
{Lie} algebras has been given by {Romdhani} \Ro;
the complex case has been considered
first by {Ancochea-Bermudez} and {Goze} \An;
complex decomposable algebras have been studied by
{Charles} and {Diakite} \Ch.
 The variety of structure constants of complex {Lie} algebras
has been examined for $n=4, 5, 6$ by {Kirillov} and {Neretin} \Ki.
{Grunewald} and {O'Halloran} \Gru\ have investigated the
complex, nilpotent {Lie} algebras for $n\leq 6$.

For $n=6$, all real nilpotent {Lie} algebras are classified
by {Morozov} \Mo;
solvable, non nilpotent
{Lie} algebras have been classified by {Mubarakzjanov}  \cMu;
and solvable real {Lie} algebras containing nilradicals are
classified by {Turkovski} \bTur, thus  completing
the classification of the solvable ones.
Both give reference
to the early
work of {Umlauf} \Um\
already classifying the nilpotent complex $6$-dimensional
{Lie} algebras.
{Mubarakzjanov} also classified
real {Lie} algebras up to $n=5$ in \bMu.
In \aMu\ he treats the case of real $n=4$, giving reference
to the early works of {Lie} \bLie\ for complex algebras with
$n\leq 4$ and, for the real $3$-dimensional case, to
{Bianchi} \bBi\ and later equivalent classifications of
{Lee} \Lee\ and {Vranceanu} \Vra.

The $3$-dimensional real {Lie} algebras,
are given by the so called {Bianchi} types,
classified independently
first by S. {Lie} \Lie\ and then by L. {Bianchi}  \Bi.
The original classification of {Bianchi}  revealed
9 inequivalent types of $3$-parameter {Lie} groups $G_3$,
numbered usually by the Roman numbers $\I,\ldots,\IX$.
The types of number \VI\ and \VII\
are actually
$1$-parameter sets
of {Lie} algebras,
\VIh\ resp. \VIIh,
with $h\geq 0$ all inequivalent.
We will refer to the inequivalent
$3$-dimensional real {Lie} algebras as the {Bianchi} types.
Our choice of the parameter $h$ is according
to {Landau-Lifschitz} \Lan,
which agrees for \VIIh with {Behr}'s choice
in {\Est}.

When the isomorphism classes of {Lie} algebras for a given
real (or complex) dimension are known in a given dimension,
one can start to compare their algebraic structure systematically
and find their algebraic characteristics, i.e. the invariants.
So  {Paterea} and {Winternitz} \PaW\
determined  subalgebra structures for  real {Lie} algebras with $n\leq 4$.
{Grigore} and {Popp} \Gri\ developed a general classification
of subalgebras of {Lie} algebras with solvable ideal,
and invariants of real {Lie} algebras have been calculated
for $n\leq 5$ by {Patera, Sharp, Winternitz} and {Zassenhaus}
{\PaSWZ}.

But the algebraic properties of {Lie} algebras are also related to
the topological structure of the space of all {Lie} algebras
in a given dimension.
On the space of all structure constants of real {Lie} algebras
in $n$ dimensions
{Segal} has introduced
(see page 255 in \Se) the subspace topology induced from $\R^{n^3}$.
The space $K^n$ of all isomorphism classes of real
$n$-dimensional {Lie} algebras
under general linear isomorphisms ${\GL}(n)$ of their generators
has a natural weakly separating
(i.e. $T_0$, not $T_1$) non-{Hausdorff}  topology $\kappa^n$,
induced as the quotient topology from the {Segal} topology
by the equivalence relation given on the structure constants
via the action of $\GL(n)$.
This topology has been discovered
and described explicitly by {Schmidt} for $n\leq 3$ in \aSch\ and
more generally in \bSch.

As a real vector space, a Lie algebra admits also a natural orientation.
By the exponential map, for any {Lie} algebra there exists an
associated
{Lie} group which similarly admits the corresponding orientation
as a differentiable manifold.

Note that throughout the following any index or property concerning
orientation is set in brackets $()$ iff the corresponding quantity
can be considered optionally with or without reference to an orientation.

The present paper is organized by the following sections.

Sec. 2 resumes some well-known facts on {Lie} algebras
and topology needed in the sequel.

Sec. 3 describes the general construction of the topological spaces
$(K^n,{\kappa}^n)$ and $(K^n_{or},{\kappa}^n_{or})$, respectively
with and without  orientation of the {Lie} algebras
as vector spaces.

Sec. 4 shows how those solvable elements
of $K^n$ which contain all the same ideal $J_{n-1}$ can be characterized
against each other by the normalized version (NJNF) of the
{Jordan} normal
form (JNF) of a single structure matrix.
Correspondingly, an oriented normalized {Jordan} normal form
(ONJNF) for the structure constants of oriented {Lie} algebras is defined.
Hence transitions $A\to B$ between {Lie} algebras can be described by
transitions between the corresponding normal forms.

Sec. 5 resumes important general properties
(see also {Schmidt} \bSch) of the
topology ${\kappa}^n_{(or)}$ and shows up further features of orientation
duality for arbitrary dimension $n$.
The structure of $K^n_{or}$, the space of equivalence classes of oriented
{Lie} algebras, as compared to its unoriented counterpart $K^n$,
has also been described in {Rainer} \aRa.
A generalization
of {Schmidt}'s notion of atoms
is made for arbitrary subsets of
$K^n_{(or)}$.
This is applied to the case of the non-selfdual subset
$K^n_{or}\setminus K^n_{SD}$,
decomposing it for $n=3$ and $n=4$
into its connected components $K^n_+$ and $K^n_-$.

Sec. 6 is devoted to the topology of the non oriented $K^n$ for $n\leq 4$.
The topological structure for $n\leq 4$ has also been described
by {Rainer} \aRa.
The $T_0$
topology $\kappa^n$ provides for $n\geq 3$ a rich local structure of $K^n$,
which we describe for $n\leq 4$.

In Sec. 6.1 the topological structure of $K^n$ for
$n\leq 3$ is analysed by use of the NJNF. So, using a quite different
method, we reproduce the results of
{Schmidt} in \aSch\ and \bSch.

Sec. 6.2 presents the detailed analysis of
the components of
$K^4$, their possible $\kappa^4$ limits, and transitions between them.
Thereby the relation between the
different classification schemes of
{Mubarakzjanov} \aMu, {Patera, Winternitz} \PaW\ and
{Petrov} \Pe\ is clarified.

Sec. 6.3
determines the topological structure of $K^4$.
Its parametrically connected components are related in a transitive
network of $\kappa^4$ transitions.

Sec. 7 is devoted to the topology of the oriented $K^n_{or}$ for $n\leq 4$,
which is also described
in {Rainer} \bRa.

In Sec. 7.1 the topological structure of $K^n_{or}$ for
$n\leq 3$ is analysed by use of the ONJNF, in correspondence with
results listed by {Schmidt} \bSch.
We give the connected components
$K^3_\pm$ explicitly.

Using the same method, Sec. 7.2 examines the orientation duality structure
of $K^4_{or}$ in detail. In particular,
we determine the connected components
$K^4_\pm$.

In Sec. 8 we discuss the present results.

\section{\bf Preliminaries}
\setcounter{equation}{0}
In the following we remind
shortly some of the notions needed throughout this paper.
A (finite-dimensional) {Lie} algebra is a
(finite-dimensional) vector space $V$, equipped with a
skew symmetric bilinear product
$[\cdot,\cdot]$ called {Lie} bracket,
which maps $(X,Y)\in V\times V$ to $[X,Y]=-[Y,X]\in V$ and satisfies
$\sum_{cycl.\atop X,Y,Z} [[X,Y],Z]=0, \forall X,Y,Z \in V$. The dimension of
the {Lie} algebra is the dimension of the underlying vector space.
Here and in the following all {Lie} algebras and vector spaces
are assumed to be finite-dimensional.
If the vector space is real resp.
complex, we say that the {Lie} algebra is real resp. complex.
If nothing else is specified in the following a {Lie} algebra or a
vector space is assumed to be real.

For a {Lie}
algebra $A$
the descending central series
of ideals is defined recursively by
\begin{equation}
C^0A:=A \quad{\rm and}\quad C^{i+1} A:=[A,C^{i}A]\subseteq C^{i}A.
\end{equation}

$A$ is called {\em nilpotent},
iff there exists a $p\in \N$, such that $C^pA=0$, i.e.
the descending central series
of ideals terminates at the zero ideal.

Furthermore for a {Lie} algebra $A$
the derivative series of ideals is defined
recursively by
\begin{equation}
A^{(0)}:=A \quad{\rm and}\quad A^{(i+1)}:=[A^{(i)},A^{(i)}]\subseteq A^{(i)}.
\end{equation}

$A$ is called {\em solvable},
iff there exists a
finite $q\in \N$, such that $A^{(q)}=0$, i.e. the derivative
series of ideals terminates at the zero ideal.

In the following, we consider
real {Lie} algebras of fixed finite dimension
$n\geq 2$ (for $n=1$ there is only 1 type of {Lie} algebra, namely
the {Abel}ian $A_1$), classified up to equivalence via real $\GL(n)$
transformations
of their linear generators $\{e_i\}_{i=1,\ldots,n}$,
which span an $n$-dimensional real vector space, which may in the following
be identified with $\R^n$ or the tangent space $T_xM$ at any point $x$ of an
$n$-dimensional smooth real manifold $M$. The {Lie} bracket $[\ ,\ ]$ is
given by its action on the generators $e_i$,
which is encoded in the structure constants $C^k_{ij}$,
\begin{equation}
[e_i,e_j]=C^k_{ij} e_k.
\end{equation}
(The sum convention is always understood implicitly, unless stated
otherwise.)
The bracket $[\ ,\ ]$ defines a {Lie} algebra, iff the structure constants
satisfy the $n\{{n\choose 2}+{n\choose 1}\}$ antisymmetry
conditions
\begin{equation}
C^k_{[ij]}=0,
\end{equation}
and the $n\cdot{n\choose 3}$ quadratic compatibility constraints
\begin{equation}
C^l_{[ij}C^m_{k]l}=0
\end{equation}
with nondegenerate antisymmetric indices $i,j,k$.
Here $_{[\quad ]}$ denotes
antisymmetrization w.r.t. the indices included.

Note that Eq. (2.5) is satisfied automatically by Eq. (2.4), if
the bracket is derived via $[e_i,e_j]\equiv e_i\cdot e_j-e_j\cdot e_i$
from an associative multiplication $e_i\cdot e_j$. In this case
Eq. (2.5) is an {\em identity}, called {Jacobi} identity.
Otherwise Eq. (2.5) is an {\em axiom}, which might be called
{Jacobi} axiom.
If there is an (adjoint) matrix representation of the algebra,
it is associative and hence
satisfies the {Jacobi} axiom (2.5) trivially, i.e. as identity.
We will not assume the existence
of any matrix representation nor any associative algebra multiplication,
because we want all the data for a {Lie} algebra to be encoded in the
structure constants. Hence we take (2.5) as an axiom.

The space of all sets $\{C^k_{ij}\}$ satisfying the {Lie} algebra
conditions (2.4) and (2.5)
can be viewed as a subvariety  $W^n \subset \R^{n^3}$ of dimension
\begin{equation}
\dim W^n \leq n^3 - \frac{n^2(n+1)}{2}  =\frac{n^2(n-1)}{2}.
\end{equation}
For $n=3$ the structure constants can be written as
\begin{equation}
C^k_{ij}=\eps_{ijl}(n^{lk}+\eps^{lkm}a_m),
\end{equation}
where $n^{ij}$ is symmetric and $\eps_{ijk}=\eps^{ijk}$ totally
antisymmetric with $\eps_{123}=1$.
With Eq. (2.7) the constraints
Eq. (2.5) are equivalent to
\begin{equation}
n^{lm}a_m=0,
\end{equation}
which are $3$ independent relations.
Actually, {Behr} has first classified the {Lie} algebras in $K^3$
according to their possible inequivalent eigenvalues of $n^{lm}$ and
values of $a_m$
(see {Landau-Lifschitz} \Lan).

With Eq. (2.8) also Eq. (2.5) is nontrivial for $n=3$.
Therefore the inequality in Eq. (2.6) is strict for $n\geq 3$.

Throughout the following, we will need the {\em separation axioms}
from topology (for further reference see also {Rinow} \Ri).
A given topology on a space $X$ is {\em separating} with increasing
strength if it satisfies one or more of the following axioms.
\nl
{\bf Axiom $T_0$}: For each pair of different points there is an open set
containing only one of both.
\hfill\mbox\break
\nl
{\bf Axiom $T_1$}: Each pair of different points has a pair of open
neighbourhoods with their intersection containing none of both points.
\hfill\mbox\break
\nl
{\bf Axiom $T_2$} ({Hausdorff}):
Each pair of different points has a pair of disjoint neighbourhoods.
\hfill$\Box$\break
It holds: $T_2\imp T_1\imp T_0$. If a topology is only $T_0$, but not $T_1$,
we say that it is only {\em weakly separating}  and speak also shortly
of the {\em weak} topology.
(The present notion {\em weak} should not be confused with another
one from functional analysis, which is not meant here.
{\em Separability} of the topological space is defined here by the
separation (german: Trennung) axioms $T_0, T_1$ or $T_2$.
This should not be confused with a further notion related to
the existence of a countable dense subset.)
The separation axioms can equivalently be characterized in terms
of sequences and their limits.

{\bf Lemma.} For a topological space $X$ the following equivalences hold:
\nl
a) $X$ is $T_0$ $\Iff$ For each pair of points there is a sequence converging
only to one of them.
\nl
b) $X$ is $T_1$ $\Iff$ Each constant sequence has at most one limit.
\nl
c) $X$ is $T_2$ ({Hausdorff}) $\Iff$ Each {Moore-Smith}-sequence has
at most one limit.
\nl
\hfill\mbox\break
(As a generalization of an ordinary sequence, a {Moore-Smith} sequence
is a sequence indexed by a (directed) partially ordered set.)
$T_1$ is equivalent to the requirement that each one-point set is closed.

We define for the following
the real {\em dimension} of a set as the largest number $k$
such that a subset homeomorphic to $\R^k$ exists.

\section{\bf Spaces $K^n$ and $K^n_{or}$ of {Lie} algebras}
\setcounter{equation}{0}
The space of structure constants $W^n$ can also be considered as a subvariety
of the fibrespace of
the tensor bundle $\wedge^2T^*M\otimes TM$
over any point of some smooth $\GL(n)$-manifold $M$.
If $M$ is oriented, the structure group of its tangent vector bundle
$TM$ is reduced from $\GL(n)$ to its normal subgroup
\begin{equation}
\GL^+(n)=\{A\in \GL(n):\det A > 0\}.
\end{equation}
Then $W^n$ gets an additional structure induced from $\wedge^2T^*M\otimes TM$
by the orientation of $M$. \nl
$\GL(n)$ basis transformations  induce
$\GL(n)$ tensor transformations between equivalent structure constants.
\begin{equation}
C^k_{ij} \sim (A^{-1})^k_h\ C^h_{fg}\ A^f_i\ A^g_j \ \ \forall A \in \GL(n),
\end{equation}
where $\sim$ denotes the equivalence relation.
This induces the space
\begin{equation}
K^n=W^n/\GL(n)
\end{equation}
of equivalence classes  w.r.t. the
nonlinear action of $\GL(n)$ on $W^n$. The analogous space for the oriented
case is
\begin{equation}
K^n_{or}=W^n/\GL^+(n).
\end{equation}
The ${\GL}(n)$ action on $W^n$ is not free in general. It holds:
\begin{equation}
\dim W^n> \dim K^n_{(or)}\geq \dim W^n - n^2.
\end{equation}
The first inequality in Eq. (3.5) is a strict one,
because the (positive) multiples of the unit matrix in $\GL^{(+)}(n)$ give
rise to equivalent points of $K^n_{(or)}$.
Eqs. (2.6) and (3.5) provide only insufficient information on
$\dim K^n$. The latter is still unknown for general $n$. (For the
analogous complex varieties {Neretin} \Ne\ has given
an upper bound estimate.)

Let $\phi_{(or)}: W^n\to K^n_{(or)}$ be the canonical map for the
equivalence relation $\sim$ defined by the action of $\GL^{(+)}(n)$
in $W^n$.
The natural topology $\kappa^n_{(or)}$ of $K^n_{(or)}$ is given as the
quotient topology of the induced subspace topology of $W^n \subset \R^{n^3}$
w.r.t. the $\GL^{(+)}(n)$ equivalence relation.

In the oriented case, orientation reversal of the basis yields a
natural $Z_2$-action on $K^n_{or}$. This action is not free in general.
Hence the fibres of the projection
\begin{equation}
\pi: K^n_{or}\to K^n = W^n/{\GL}(n)=K^n_{or}/Z_2
\end{equation}
can be either $Z_2$ or $E$. In the first case
there is a pair of dual points, i.e. points that transform into each
other under the $Z_2$-action,
in the latter case it is a selfdual point in $K^n_{or}$.
The latter therefore
decomposes into a selfdual part $K^n_{SD}$, on  which $Z_2$ acts
trivially, and 2 conjugate parts $K^n_{\pm}$. The latter are isomorphic to
each other by that reflection in $GL(n)$ that is chosen to define $Z_2$ in
Eq. (3.6).
\begin{equation}
K^n_{or}=K^n_{SD}\oplus K^n_{+} \oplus K^n_{-},
\end{equation}
where $\oplus$ denotes the disjoint union of subvarieties.

The projection $\pi$ has
the property that its restriction to $K^n_{SD}$ is the identity.
Therefore it is useful to make the following
\nl

{\bf Definition 1.}
A point $A\in K^n$ is called {\em selfdual} if
$\pi^{-1}(A)\subset K^n_{or}$ consists of a single point, and
{\em non-selfdual}
if $\pi^{-1}(A)$ consists of a pair of dual points,
denoted by $A^R$ and $A^L$
respectively.
\hfill$\Box$\break
In order to yield a more explicit notion of selfduality, we formulate
\nl

{\bf Lemma.}
{\em
A {Lie} algebra $A$ is selfdual, $A\in K^n_{SD}$,
if and only if there exist two different bases of $\R^n$ possessing different
orientation such that all the structure constants $C^k_{ij}$ concerning both
bases coincide.
}
\hfill$\Box$\break

Obviously a direct sum of a selfdual algebra with any other algebra
is selfdual.

Let us mention already here that
$K^n_{SD}$ is nonvoid for $n\geq 1$ while
$K^n_\pm$ are nonvoid sets only for $n\geq 3$.
We will see in Sec. 5 and 7 that the latter are actually nonvoid
for $n=3,4,5$ and at least any further odd $n$.
In any case $K^n_\pm$ are connected to $K^n_{SD}$.
We will see  in Sec. 7 that each of $K^n_\pm$ is connected
for $n=3$ and $n=4$.

Note that for each pair of conjugate {Lie} algebras $A^R$ and $A^L$ it
is a priori completely arbitrary which one is assigned to $K^n_+$ and which
one to $K^n_-$. In order to reduce this arbitrariness, in Sec. 4 we will
minimize the number of connected components of $K^n_\pm$ to a single
component each, thus making $K^n_+$ and $K^n_-$ disconnected to each other.
However this requires first a better understanding of the
topological structure $K^n_{or}$.

When we do not want to care about effects of orientation,
instead of {Schmidt}'s topological space
$(G_n,\tau)\equiv (K^n_{or},\kappa^n_{or})$ from \bSch\
we will consider its projection  to $(K^n,\kappa^n)$ by Eq. (3.6).

Let us define now the notion of transitions $A\to B$ in
$K^n_{(or)}$.
\nl

{\bf Definition 2.}
Consider $A, B \in K^n_{(or)}$ with $A\neq B$.
If there is a sequence $\{A_i\}_{i\in \N}$ with
$A_i=A$ for all $i\in \N$  which for $i\to \infty$
converges to $B$ in the topology $\kappa^n_{(or)}$,
we say that there is a {\em transition} $A\to B$ in the topology
$\kappa^n_{(or)}$.
\hfill$\Box$ \break
Note that this definition makes sense because $K^n$ is a $T_0$
but not a $T_1$ space.
A transition is a special kind of limit characteristic for this topology.
\nl

{\bf Convention.}
We distinguish in notation
between a concrete realization of a {Lie} algebra, $A$, and its
equivalence class, $[A]$, where ever this is relevant.
In the following, the former will an adjoint representation
of the latter, sometimes also called abstract,
{Lie} algebra.
However for notational simplicity we prefer to denote a point in $K^n$
by $A$ rather than by $[A]$. If the context does not give
the opportunity for confusion, $A$ is implicitly understood as
a shorthand for the (abstract) {Lie} algebra $[A]$.
\hfill$\Box$ \break
In the topology $\kappa^n_{(or)}$, a transition $A\to B$
occurs if and only if $B\in \cl \{A\}$.
For this transition the source $A$ is not closed, and the target
$B$ is not open in any subset of $K^n_{(or)}$ containing both of them.
In general, a point of $K^n_{(or)}$ will be  neither open nor closed.
Open points  only appear as a source, and not as a target, of transitions.
The structure of the rigid {Lie} algebras, which correspond just
to these open points, is examined in {Charles} \Cha.

Special kinds of transitions
on a certain 2-point set $\{A,B\}$
of {Lie} algebra isomorphism classes
are the contractions of {In\"on\"u-Wigner} \In\ and their
generalization by {Saletan} \Sal. For convenience let us define
these here.

Consider a $1$-parameter set of matrices $A_t\in\GL(n)$ with
$0<t\leq 1$, having a well defined matrix limit $A_0:=\lim_{t\to 0} A_t$
which is singular, i.e. $\det A_0=0$.
For given structure constants $C^k_{ij}$ of a {Lie} algebra
$A$ let us define for
$0<t\leq 1$ further structure constants
$C^k_{ij}(t):=(A^{-1}_t)^k_h\ C^h_{fg}\ (A_t)^f_i\ (A_t)^g_j$, which
according to (3.2) all describe the same {Lie} algebra $A$.
If there is a well defined limit $C^k_{ij}(0):=\lim_{t\to 0} C^k_{ij}(t)$
satisfying conditions (2.4) and (2.5)
then this limit defines structure constants of a {Lie} algebra $B$,
and the associated limit of {Lie} algebras $A\to B$
is called {\em contraction} according to {Saletan} \Sal\
or briefly {Saletan} {\em contraction}.
Note that a {Saletan} contraction $A\to B$ might yield either
$B=A$, then it is called {\em improper}, or $B\neq A$, then it is
a transition of {Lie} algebras.

A {Saletan} contraction is called
{In\"on\"u-Wigner} {\em contraction}
if there is a basis $\{e_i\}$ in which
$$
A(t)=
\left(
\begin{array}{cc}
 E_m & 0 \\
 0   & t\cdot E_{n-m}
\end{array}
\right)
\qquad \forall t\in [0,1],
$$
where $E_k$ denotes the $k$-dimensional unit matrix.
This definition closely follows {Conatser} \Co.

Given the latter decomposition, {In\"on\"u} and {Wigner} \In\ have
shown that the limit $C^k_{ij}(0)$ exists iff $e_i, i=1,\ldots,m$
span a subalgebra $W$ of $A$, which then characterizes the
contraction.

{Saletan} \Sal\ gives also a technical criterion for the
existence of the limit
$C^k_{ij}(0)$ defining his general contractions.

We only remark here that,
while a general {Saletan} contraction might be nontrivially iterated,
the iteration of an {In\"on\"u-Wigner} contraction is always improper,
i.e. no further contraction takes place.

Not every transition $A\to B$
corresponds to an {In\"on\"u-Wigner} contraction.
We will see some examples of transitions, which are given
only by a more general {Saletan} contraction \Sal.
However we will find also
transitions $A\to B$,
which are not even given by a {Saletan} contraction.

Transitions $A\to B$ in the topology $\kappa^n$ reveal for $n\geq 3$ a more
complicated structure of the underlying space $K^n$.
In $K^n$ transitions $A\to B$ and $B\to C$ imply a transition $A\to C$;
this means that transitions are transitive.
There is a partial order,
$A \geq B :\Iff B\in \cl \{A\} \Iff A\to B$
(which is also called the {\em specialization order}), which gives
$K^n_{(or)}$ the structure of a transitive network of transitions.
Since {Saletan} contractions \Sal\ are not transitive they do not
exhaust all kinds of possible $\kappa^n$ transitions.

Given the topology of $K^n$, on any 2-point subset $\{X,Y\}\subset K^n$
we can take the induced topology and consider the set
$T^n:=\{\{X,Y\}\subset K^n\vert X\neq Y\}$ of all 2-point topological
subspaces of $K^n$.
Note that a $T_0$ topological space, like that of $K^n$ for $n\geq 3$,
is in general not determined by the set $T^n$ of all its induced
2-point topological subspaces.
However if the topological space under consideration
is finite then $T^n$ determines already its topology,
which is trivially true
for $K^1$ and $K^2$.

\section{\bf Normal forms of the structure constants}
\setcounter{equation}{0}
The structure constants of $A_n\in [A_n]\in K^n$ are given by
the $n$ matrices $C_i:=(C^k_{ij})$, $i=1,\ldots,n$,
with rows $k=1,\ldots,n$ and
columns $j=1,\ldots,n$. $C_i$ is just the matrix of ad$e_i$ w.r.t.
the basis $e_1,\ldots,e_n$.\nl
By Eq. (2.4), the column $j=i$ vanishes identically $\forall C_{i}$.
Furthermore
the diagonals $(C^j_{ij})$, $i=1,\ldots,n$ (no j-summation),
determine the rows
with $k=i$, since $(C^i_{ij})=(-C^i_{ji})$, $i=1,\ldots,n$
(no i-summation).
Therefore $A_n$ is described completely by the
$(n-1)\times (n-1)$-matrices
$C_{<i>}:=(C^k_{ij})$, $i=1,\ldots,n$, with $k,j\neq i$ and
$1\leq k,j\leq n$.\nl
In the special case where $A_n$ has an ideal
$J_{n-1}\in [J_{n-1}]\in K^{n-1}$,
we take without restriction $[A_n]/[J_{n-1}]=\mbox{span}(e_n)$.
Then $A_n$ with a given
$J_{n-1}$ is described completely by $C_{n}$ or $C_{<n>}$ only.
\nl

{\bf Definition 3.} The {\em normalized} JNF (NJNF) of a matrix $C$ is given
by the {Jordan} normal form, abbreviated JNF,
of $C$ modulo $\R\setminus \{0\}$, i.e. given by the equivalence
class of JNFs, which differ only by a common absolute scale and a common
overall sign of their nonzero eigenvalues w.r.t. the eigenvalues of $C$.
(The {Jordan} block structure and
the multiplicities are the same for all of
them.)
\hfill $\Box$\break
Thus a normalization convention for the JNF is the division of all
nonzero eigenvalues by a fixed element of $\R\sm\{0\}$. If not stated
otherwise,
we divide in the following just by the (absolutely) largest eigenvalue in order
to
represent the NJNF class of the JNF.\nl
Note that the $n^{th}$ row and column of $C_n$
add only an additional
eigenvalue 0 (as {Jordan} block) to the JNF or NJNF of $C_{<n>}$.
Since absolute scaling of all eigenvalues of a structure matrix $C_{<n>}$
by $\lambda \in \R\sm\{0\}$ can be
achieved by stretching the basis $\{e_i\}$ homogeneously by $\lambda^{-1}$,
it is an equivalence transformation of the algebra. On the other hand it is
evident that changing in $C_{<n>}$ the ratio $r$ of any 2 eigenvalues to
$r'$, such that $r'$ is not a ratio of any original
eigenvalues, changes the equivalence class.

{\bf Theorem.}
{\em
Consider the set of algebras $A_n$ which have a
common (abstract) ideal $J_{n-1}$.
Then $A^{(1)}_n\sim A^{(2)}_n$, iff the matrices $C^{(1)}_{<n>}$ and
$C^{(2)}_{<n>}$ have the same \NJNF.
}\nl
{Proof:} $A^{(1)}_n\sim A^{(2)}_n$ iff $\exists M\in {\GL}(n):
{C^{(1)}}^k_{ij} = (M^{-1})^k_h\ {C^{(2)}}^h_{fg}\ M^f_i\ M^g_j
\sim M^f_i {C^{(2)}}^h_{fg}$. By linearity of $[\ ,\ ]$ in the second
argument, the linearly independent recombinations
${\tilde{C}}^{(2)}_i:=M^f_i\ C^{(2)}_{f}$ describe still the same algebra
as $C^{(2)}_{i}$. In particular, the (abstract) ideal $J_{n-1}$
is invariant under $M$.
Since the algebras have the same ideal $J_{n-1}$,
they are characterized by
the matrices $C^{(1)}_{<n>}$ resp. $C^{(2)}_{<n>}$. They describe
inequivalent algebras, iff $C^{(1)}_{<n>}$ is inequivalent (modulo overall
scaling by $M=\lambda E_n,\ \lambda \in \R\sm\{0\}$) to
${\tilde{C}}^{(2)}_{<n>}$ and therefore also to $C^{(2)}_{<n>}$.
But the equivalence class of any
structure matrix $C_{<n>}$ is described by its (real) JNF modulo homogeneous
scaling of the eigenvalues with $\lambda \in \R\sm\{0\}$.
\hfill $\Box$ \break
\nl
Already {Mubarakzyanov} \aMu\ had realized the advantage given by
an  ideal $J_{n-1}$ of codimension $1$. Since then also others, like
{Magnin} \Mag\ within the nilpotent {Lie} algebras of dimension $\leq 7$,
systematically cosidered subclasses of algebras which have a fixed
{Lie} algebra of codimension $1$.

In the following, we consider without restriction
of generality the ideals $J_{(n-1)}$
in the normal form given by the NJNF of the structure constants.
The equivalence class $[A_{n}]$ of any algebra $A_{n}$ with a
fixed normal class ideal and additional structure constants from $C_{<n>}$
will be characterized in the following by the NJNF of $C_{<n>}$
and denoted by
$$
\NJNF(A_n) := \NJNF(C_{<n>}).
$$
Now we can define the ONJNF of structure matrices $C_{<n>}$ of  oriented
{Lie} algebras.
\nl

{\bf Definition 4.}
The ONJNF of the structure matrix $C_{<n>}$ of an oriented {Lie} algebra
$A_n$ is set identical to its NJNF if $A_n$ is selfdual, and
it is given as
$\ONJNF(C_{<n>}):=\pm \NJNF(C_{<n>})$ for $A_n\in K^n_{\pm}$ respectively.
\hfill $\Box$\break
If $A_n$ is characterized by an ideal $J_{n-1}$ in normal form
then we set
$$
\ONJNF(A_n) := \ONJNF(C_{<n>}).
$$

\section{\bf General properties of $\kappa^n$ and $\kappa^n_{or}$}
\setcounter{equation}{0}
In this section we describe the general topological properties of the
topological space $(K^n_{(or)},\kappa^n_{(or)})$.
Let us first remind some general properties
from {Schmidt} \bSch\ (where also more details and proofs
can be found).
\nl

{\bf Proposition.}
{\em
$K^n_{(or)}$ has the following properties w.r.t. $\kappa^n_{(or)}$:
}
\nl
a)
{\em
The {Abel}ian algebra $\{nA_1\} \subset K^n_{(or)}$ is the only closed
1-point set and is contained in any nonempty closed subset of $K^n_{(or)}$.
}
\nl
b)
{\em
$K^n_{(or)}$ is connected and compact.
}
\nl
c)
{\em
$K^n_{*(,or)}:=K^n_{(or)}\sm\{nA_1\}$ is a compact space, but
$K^n_{*(,or)}$ is not a
closed subset of $K^n_{(or)}$.
}
\nl
d)
{\em
$K^n_{*(,or)}$ is {Hausdorff} $(T_2)$ for $n=2$ only.
}
\nl
e)
{\em
For $n\geq 2$ (resp. $n\geq 3$) the separability of $K^n_{(or)}$
(resp. $K^n_{*(,or)}$) is only weak ($T_0$, i.e. for each pair of points
there is a sequence converging to only one of them).
}
\hfill$\Box$\break
%
d) and e) correspond to the fact that, though $K^n$ is still an algebraic
variety (defined by purely algebraic relations (2.4), (2.5) and
(3.2)), it can not be expected to be a (topological $T_1$)
manifold.
$K^n$ is the orbit space of $W^n$ w.r.t. the action of the
noncompact group ${\GL}(n)$, which behaves algebraically badly on $W^n$
for $n\geq 2$. So some of the orbits (the elements
of $K^n$) are closed in $K^n$, others are not.

Strong separability ($T_1$, i.e.  each constant sequence has at most one
limit) would imply that there should not exist
transitions $A \to B$ between
inequivalent {Lie} algebra classes $A\not\sim B$, given by a sequence
$\{A_i\}$ of {Lie} algebras of class $A$ converging to a
{Lie} algebra of class $B$.
But this is exactly what happens for dimension $n\geq 2$, as will be
seen explicitly below. Obviously transitions $A\to B$ will be transitive,
which decisively effects the topology of $K^n$.

\nl
Transitions which are impossible in a given dimension $n$
can become possible
after {Abel}ian embedding into dimension $n+1$. Therefore the following
lemma holds.\nl

{\bf Lemma 1.}
{\em
The {Abel}ian embedding $\oplus \R$ of $K_n$
into $K_{n+1}$ is continuous, but for $n\geq 2$ not homeomorphic.
}
\hfill$\Box$ \break

So we are led to the following
\nl

{\bf Definition 5.}
The {\em essential dimension}  of an $n$-dimensional (oriented)
{Lie} algebra $A_n$ is
defined as the smallest possible number $n_e\leq n$, such that
$A_n=A_{n_e}\oplus \R^{n-n_e}$. The essential-dimensional
subset of $K^n$ is defined as $K^n_{de}=\{A\in K^n\vert n_e(A)=n\}$
\hfill$\Box$ \break

{\bf Lemma 2.}
{\em
The subsets $\{A\in K^n\vert n_e(A)\leq m\}$ for any
fixed $m\leq n$ need not to be closed.
}
\hfill$\Box$ \break
This is due to the existence of transitions or limits
of structure constants in NJNF
such that one or more NJNF eigenvalues degenerate to another one
(in Lemma 2 it is the eigenvalue $0$), initially distinct from them;
in this case the algebraic multiplicity of this eigenvalue increases
automatically, but its geometric multiplicity (expressed by its number of
{Jordan} blocks) may remain constant, since an eigenvector of
an eigenvalue different from the limit eigenvalue may converge to a
principal (not necessarily eigen) vector of the limit eigenvalue
($0$ for Lemma 2).\nl

A {Lie} algebra characterized by structure constants $C^k_{ij}$ is
called {\em unimodular}
(on a corresponding {Lie} group) iff $\trace(C_i) = C^k_{ik} =0\
\forall i$, where the adjoint representation is generated by the matrices
$C_i$. We denote the subset of all points in $K^n$ that correspond
to unimodular {Lie} algebras by $U^n$, and set $U^n_*=U^n\cap K^n_*$.
Since the zero set of a continuous function is always closed, we have\nl

{\bf Lemma 3.}
{\em
The unimodular subset $U^n\subset K^n$ is closed and compact.
}
\hfill$\Box$ \break

For $n\geq 2$ the structure constants of any {Lie} algebra admit,
as a tensor $C$, the irreducible decomposition \bSch
\begin{equation}
C^k_{ij}=D^k_{ij}+\delta^k_{[i}v_{j]}
\end{equation}
in a tracefree part $D$ with tensor components $D^k_{ij}$ (the trace free
condition for $D$ can be written as
$\trace(D_i)=D^k_{ik}=0\ \forall i$),
and a vector part, constructed from a vector $v$ with components
$v_i=C^j_{ij}/(1-n)$ and the Kronecker
symbol of components $\delta^k_i$ (remind the sum convention over upper and
lower indices and the convention to perform an antisymmetric sum over all
permutations of the indices included in $_{[\quad ]}$).
The {Lie} algebra is {\em unimodular}
(like any associated connected {Lie} group),
iff it is tracefree, $v\equiv 0$, and it is said to be of
{\em pure vector type},
iff $D\equiv 0$. In this sense the unimodular and pure vector type
are complementary.

The class ${\V}^{(n)}$ of pure vector type is selfdual for all $n\geq 2$.
It is the generalization of the unique non-{Abel}ian $2$-dimensional
algebra $A_2$ (see Sec. 6) to arbitrary $n$.
So for each $n$, there exists exactly
one non-{Abel}ian pure vector type {Lie} algebra, denoted by $\V^{(n)}$ because
for $n=3$ it is the Bianchi type \V.
It has the {Abel}ian ideal $\I^{(n-1)}$ and
$[\NJNF(\V^{(n)})]^k_j=\delta^k_j$.
For convenience we mention explicitly the nonvanishing commutators
of $\V^{(n)}$, for an adapted basis $\{e_1,\ldots,e_n\}$:
\begin{equation}
[e_n, e_i] =e_i,\ i=1,\ldots,{n-1}.
\end{equation}

The $3$-dimensional Heisenberg algebra (= Bianchi type \II)
is defined in its NJNF by the
nonvanishing commutators
\begin{equation}
[e_3, e_2] =e_1.
\end{equation}
By {Abel}ian embedding we define the
class ${\II}^{(n)}:={\II}\oplus \R^{n-3}$ for $n\geq 3$.
Like\II, it is unimodular and nilpotent of degree 2.
${\II}^{(n)}$ is non-selfdual for $n=3$ and selfdual for $n\geq 4$.
Its nonvanishing structure constants for an adapted basis
$\{e_1,\ldots,e_n\}$)
are given by Eq. (5.3)
with all indices increased by $n-3$.


In {Schmidt} \bSch, an element $A_n\in K^n_*$ for which its closure
in $K^n$ consists of 2 elements only, $\cl\{A_n\}=\{A_n,nA_1\}$,
was called an {\em atom}.
Here we will prefer to call equivalently $A_n$ an atom of $K^n_*$, iff
its closure in $K^n_*$ is $\cl_{K^n_*}\{A_n\}=\{A_n\}$.
Let us generalize this:
\nl

{\bf Definition 6.}
For any subset $S\subset K^n_{(or)}$,
an element $A\in S$ is called an {\em $S$-atom},
iff it is closed w.r.t. $S$, i.e.
$\cl_S\{A\}=\{A\}$.
\hfill$\Box$ \break
In the following, we call an $S$-atom also synonymously
an {\em atom of $S$} and assume $S=K^n_*$ if not specified otherwise.
Recall
from {Schmidt} \bSch
\nl

{\bf Theorem 1.}
{\em
For $n=2$ there is only 1 atom, $A_2\equiv {\V}^{(2)}$.\nl
For each $n\geq 3$ there exist exactly 2 atoms, the unimodular ${\II}^{(n)}$
and the pure vector type ${\V}^{(n)}$.
}
\hfill$\Box$ \break
For $n\neq 3$ all atoms are selfdual.
If we consider the corresponding atoms of $K^n_{*,or}$, then only for
$n=3$ there is a difference to the nonoriented case. Instead of the
unique non-selfdual atom $\II$ in $K^3$, there exist 2 non-selfdual
atoms, $\II^R$ and $\II^L$, in $K^3_{or}$.

For each $n\geq 3$ there is an algebra ${\IV}^{(n)}$,
given by
$[\NJNF(\IV^{(n)})]^k_j=\delta^k_j+\delta^k_{n-2}\delta^{n-1}_j$ w.r.t. to
the {Abel}ian ideal $\I^{(n-1)}$. It is selfdual for $n\geq 4$ and
non-selfdual for $n=3$.
For convenience we mention explicitly the nonvanishing commutators
of $\IV^{(n)}$, for an adapted basis $\{e_1,\ldots,e_n\}$ given by
\begin{equation}
[e_n, e_i] =e_i,\ i=1,\ldots,n-2, \quad [e_n, e_{n-1}] =e_{n-2}+e_{n-1}.
\end{equation}

$K^n_*$ is generated by infinitesimal deformations of the atoms;
this means: $K^n_*$ itself is the only open subset of $K^n_*$
which contains all atoms.
Since both, ${\IV}^{(n)} \to {\II}^{(n)}$ and
${\IV}^{(n)} \to {\V}^{(n)}$, it follows that $K^n_*$ is connected. \nl
Remark: Connectedness is trivial for $K^n$, but non-trivial for $K^n_*$.

To understand better where the exceptionality of $n=3$ w.r.t. to duality
comes from,
realize that $n_e(\V^{(n)})=n$ but $n_e(\II^{(n)})=3$
for all $n\geq 3$. In particular, $\II^{(n)}$ has essential dimension
$n_e=n$ only for $n=3$; for $n\geq 4$ it is decomposable and hence
selfdual.

More generally there holds
\nl

{\bf Lemma 4.}
{\em
$K^n_{NSD}:=K^n\sm K^n_{SD}$ is contained in the subset $K^n_{de}$
of $K^n$ for which $n_e=n$.
}
\hfill$\Box$ \break
To overcome the difference in the essential dimension of the atoms for
$n\geq 4$, let us search for atoms w.r.t. the subset $K^n_{de}$
of essential dimension $n_e=n$ in $K^n$. We find
\nl

{\bf Theorem 3.}
{\em
The set $K^n_{de}$ has the following atoms:
\nl
a)
For $n\geq 2$
exactly $1$ pure vector type atom, called $\ve(n)$.
}
\nl
b)
{\em
For $n\geq 3$ a nilpotent unimodular atom,  called $\ii(n)$,
located in the subset of algebras with ideal $\I^{(n-1)}$.
}
\nl
c)
{\em
For $n\geq 5$ further $]\frac{2}{3}(n-4)[$ mixed type atoms,
denoted $a_m(n)$, $m=2+[\frac{n-4}{3}],\ldots, n-3$,
all located in the subset of algebras with ideal $\I^{(n-1)}$.
\nl
(Here $[x]$ resp. $]x[$ denotes the largest/smallest integer
less/greater or equal than x.)

Within the subspace $K^n_{de\vert\I^{(n-1)}} \subset K^n_{de}$
given by $K^n_{de}$-algebras with ideal $\I^{(n-1)}$
there are no further $K^n_{de}$-atoms than that of {\rm a), b)} and
{\rm c)}.

$K^n_{de\vert\I^{(n-1)}}$ is connected.
}
\nl
Proof:
a) By Theorem 1 the algebra $\V^{(n)}$ is an atom of $K^n_*$.
Since $K^n_{de}\subset K^n_*$ and $n_e(\V^{(n)})=n$,
it follows
that $\V^{(n)}$ is an atom of $K^n_{de}$. Any algebra with
only nonzero components $v_i$ in the vector $v$
of the decomposition (5.1) has a transition
or limit to $\V^{(n)}$.
Hence $\ve(n):=\V^{(n)}$ is the unique (pure) vector type $K^n_{de}$-atom.

b)
Some of the algebras with some vanishing component $v_i$
have transitions or limits to an algebra with $v\equiv 0$.
Hence we have to search for unimodular $K^n_{de}$-atoms of essential
dimension $n_e=n$. Such an atom is the nilpotent algebra
$\ii(n)$ with $\NJNF(\ii(n))$
w.r.t. the ideal $I^{(n-1)}$ given for even $n$
as a direct sum
of $1$ block of $\NJNF(\ii(4))$ and further blocks of $\NJNF(\ii(3))$,
and for odd $n$
as a direct sum of $\NJNF(\ii(3))$ blocks only,
where
$$
\NJNF(\ii(3)):=
\left(
\begin{array}{cc}
 0 & 1 \\
   & 0
\end{array}
\right)
$$
and
\begin{equation}
\NJNF(\ii(4)):=
\left(
\begin{array}{ccc}
 0 & 1 &   \\
   & 0 & 1 \\
   &   & 0
\end{array}
\right).
\end{equation}
The algebra $\ii(n)$ is essential-dimensional,
because any of its subalgebras
$\ii(3)$ and $\ii(4)$ is so; it is an
$K^n_{de}$-atom, because its only possible limits
necessarily
generate a $1\times 1$-block $(0)$ in its NJNF, thus decreasing $n_e$
at least by $1$.

c)
The mixed atoms can be characterized by their NJNF w.r.t. the ideal
$\I^{(n-1)}$. Let us set
$$
\NJNF(a_m(n)):=\NJNF(\ve(m+1))\oplus\NJNF(\ii(n-m)).
$$
Since for
$m=2+[\frac{n-4}{3}],\ldots, n-3$ the geometric multiplicity $m$
(= the number of {Jordan} blocks) of the eigenvalue $1$ is
bigger than that of the eigenvalue $0$,
any transition yields an additional {Jordan} block
$0$ and hence leaves $K^n_{de}$. So, being essential-dimensional,
$a_m(n)$ is an $K^n_{de}$-atom for  $m=2+[\frac{n-4}{3}],\ldots, n-3$.

Any algebra of $K^n_{de\vert\I^{(n-1)}}$ has a combination of
transitions and parametric limits leading
to at least one of the atoms from a), b) or c), depending on the
degeneracy of its eigenvalues. The only nontrivial case,
which remains to be checked, are the algebras with
their NJNF w.r.t. an ideal $\I^{(n-1)}$ given as
$\NJNF(\ve(m+1))\oplus\NJNF(\ii(n-m))$ where
$m=1, \ldots, 1+[\frac{n-4}{3}]$ and $n\geq 4$.
But any of these has a transition to $\ii(n)$.

Let us now consider some algebra in $K^n_{de\vert\I^{(n-1)}}$ with
only nondegenerate nonzero eigenvalues. By continuous deformation
of its eigenvalues, such that every deformed algebra remains in
$K^n_{de\vert\I^{(n-1)}}$, and transitions within
$K^n_{de\vert\I^{(n-1)}}$ each of the atoms a), b) and c)
can be reached. Since these have just been seen to be the only
atoms of $K^n_{de\vert\I^{(n-1)}}$ it follows that
$K^n_{de\vert\I^{(n-1)}}$ is connected.
\hfill$\Box$ \break
$K^n_{de}$ itself might have further atoms
located in $K^n_{de}\setminus K^n_{de\vert\I^{(n-1)}}$.
Since these are difficult to find, in general
one cannot see whether $K^n_{de}$ is connected.

The nonvanishing commutators
of $\ii(n)$, $n\geq 3$,
are given for an adapted basis $\{e_1,\ldots,e_n\}$
explicitly by
$$
[e_n, e_2] =e_{1},\ [e_n, e_3] =e_{2},\ i=2,\ldots,n-1.
$$
\begin{equation}
[e_n, e_{2i+3}] =e_{2i+2},\quad  i=1,\ldots,\frac{n-4}{2}.
\end{equation}
for $n$ even and by
\begin{equation}
[e_n, e_{2i}] =e_{2i-1},\quad  i=1,\ldots,\frac{n-1}{2}.
\end{equation}
for $n$ odd.

$\ii(3)\equiv \II$ is the {Heisenberg} algebra.
The number of $\NJNF(\ii(3))$ blocks in its NJNF is even for
$n\equiv 0 \,\mbox{mod}\, 4$ or $n\equiv 1 \,\mbox{mod}\, 4$, and it is
odd for $n\equiv 2 \,\mbox{mod}\, 4$ or $n\equiv 3 \,\mbox{mod}\, 4$.

The mixed types $a_m(n)$, $n\geq 5$,
have respective algebraic and geometric multiplicities
$m=2+[\frac{n-4}{3}],\ldots, n-3$
for the eigenvalue $1$.

Their nonvanishing commutators
are given
w.r.t. an adapted basis
\nl
$\{e_1,\ldots,e_n\}$ as
$$
[e_n, e_1] =e_{1},\ldots, [e_n, e_m] =e_{m},
$$
$$
[e_n, e_{m+2}] =e_{m+1},\ [e_n, e_{m+3}] =e_{m+2},\quad i=2,\ldots,n-1,
$$
\begin{equation}
[e_n, e_{2i+m+3}] =e_{2i+m+2},\quad  i=1,\ldots,\frac{n-m-4}{2},
\end{equation}
for $n-m$ even, and by
$$
[e_n, e_1] =e_{1},\ldots, [e_n, e_m] =e_{m},
$$
\begin{equation}
[e_n, e_{2i+m}] =e_{2i+m-1},\quad  i=1,\ldots,\frac{n-m-1}{2},
\end{equation}
for $n-m$ odd.

The reflection $e_1\to -e_1$ leaves $\ve(n)$ and any mixed type
atom $a_m(n)$ invariant; hence all these atoms are selfdual.
The nilpotent
atom $\ii(n)$
remains as the only possibility for a non-selfdual $K^n_{de}$-atom
within $K^n_{de\vert\I^{(n-1)}}$.
Therefore, next we want to examine the orientation duality of  $\ii(n)$.
\nl

{\bf Theorem 4.}
{\em
For $n\geq 3$ the $K^n_{de}$-atom $\ii(n)$ is non-selfdual
only if $n\equiv 3\, \mbox{\rm mod}\, 4$.
\nl
$\ii(n)$ non-selfdual for $n\equiv 3\, \mbox{\rm mod}\, 4$
implies that $\ii(n)$ is a $K^n_{de}$-atom.
}\nl
{Proof:} A combination of the reflections $e_n\to-e_n$
and $e_{2i}\to-e_{2i}$ for $i=1,\ldots,[\frac{n-1}{2}]$
leaves $\ii(n)$ invariant.
The total number of these reflections is $[\frac{n+1}{2}]$,
which is odd for $n\equiv 1 \,\mbox{mod}\, 4$
or $n\equiv 2 \,\mbox{mod}\, 4$.
Furthermore for $n$ even, $e_i\to-e_i$, $i=1,\ldots,n-1$ yields
a reflection keeping $\ii(n)$ invariant.
So for all $n$ but $n\equiv 3 \,\mbox{mod}\, 4$ the algebra is selfdual.

Any limit of $\ii(n)$ is selfdual, because it is a $K^n_{de}$-atom
and any non-essential-dimensional algebra is decomposable and hence
selfdual.
Therefore non-selfduality for $n\equiv 3\, \mbox{\rm mod}\, 4$
implies that $\ii(n)$ is a $K^n_{de}$-atom.
\hfill$\Box$ \break
For $n\equiv 3 \,\mbox{mod}\, 4$ it was impossible to construct a reflection
leaving $\ii(n)$ invariant. But when there is no such reflection
the algebra is non-selfdual. 

Let us define for $n\geq 3$ an algebra $\iv(n)$
given for an adapted basis  $\{e_1,\ldots,e_n\}$ by
the nonvanishing commutators
$$
[e_n, e_1] =e_{1},\
[e_n, e_2] =e_{1}+e_{2},\ [e_n, e_3] =e_{2}+e_{3},\ i=2,\ldots,n-1,
$$
\begin{equation}
[e_n, e_{2i+3}] =e_{2i+2}+e_{2i+3},\ i=1,\ldots,\frac{n-4}{2},
\end{equation}
for $n$ even, and by
\begin{equation}
[e_n, e_{2i}] =e_{2i-1}+e_{2i},\ i=1,\ldots,\frac{n-1}{2},
\end{equation}
for $n$ odd.

By similar considerations as for $\ii(n)$ in Theorem $3$ one finds
that $\iv(n)$ is non-selfdual for only for $n$ odd and selfdual for
$n$ even.
In any case it has an ideal $\I^{(n-1)}$
and for $n$ odd the geometric multiplicity of its eigenvalue
$1$ of the NJNF w.r.t. $I^{(n-1)}$ can only be increased by yielding
at least two $1\times 1$ blocks of that eigenvalue, hence the
resulting algebra of such a transition is selfdual.
Apart from limits which increase multiplicity, the only further limits
of $\iv(n)$ are transitions with the eigenvalue becoming $0$,
either to $\ii(n)$ or some limit thereof.
But, according to
Theorem 4, for $n\not\equiv\,\mbox{mod}\,4$, the algebra $\ii(n)$ is selfdual.
Any limits of $\ii(n)$ are selfdual, because it is a $K^n_{de}$-atom
and any non-essential-dimensional algebra is decomposable and hence
selfdual.
Hence non-selfduality of $\iv(n)$ for $n$ odd implies
that  $\iv(n)$ is a $K^n_{NSD}$-atom for $n\equiv 1 \,\mbox{mod}\, 4$.
Non-selfduality of $\ii(n)$ for $n\equiv 3\, \mbox{\rm mod}\, 4$
implies further that $\iv(n)$ is
no atom for $n\equiv 3\, \mbox{\rm mod}\, 4$.

If for $n\equiv 3\, \mbox{\rm mod}\, 4$ resp. $n$ odd the algebras
$\ii(n)$ resp. $\iv(n)$ are in fact non-selfdual,
we get the
\nl

{\bf Corollary.}
{\em
For odd $n\geq 3$ the set $K^n_{NSD}$ has
an atom, located within
the subspace $K^n_{NSD\vert\I^{(n-1)}}$ of non-selfdual algebras
with ideal $I^{(n-1)}$.

For $n\equiv 3\, \mbox{\rm mod}\, 4$ the atom is
nilpotent unimodular, given by $\ii(n)$,
and for  $n\equiv 1\, \mbox{\rm mod}\, 4$ it is given by $\iv(n)$.
}
\hfill$\Box$ \break

The selfduality of the $K^n_{de}$-atoms $a_m(n)$ and  $\ve(n)$
excludes them as candidates
for  $K^n_{NSD}$-atoms.
It remains an open problem to determine
at least some $K^n_{NSD}$-atom for arbitrary even $n$,
and all $K^n_{NSD}$-atoms for arbitrary
$n$.
For odd $n$, besides $\ii(n)$ or $\iv(n)$, there might be further
$K^n_{NSD}$-atoms, even within
$K^n_{NSD\vert\I^{(n-1)}}$.

However, assume we succeed for some $n$ to determine all
$K^n_{NSD}$-atoms and
furthermore to show that $K^n_{NSD}$ is connected for that $n$.
In Sec. 7 we will actually see that, for $n=3$ the
{Heisenberg} algebra
$\ii(3)$ is the only non-selfdual atom, hence $K^3_{NSD}$ is connected,
and for $n=4$,
with the topology of $K^4$ obtained in Sec. 6.3
and the non-selfdual algebras of Sec. 7.2.2,
the resulting non-selfdual set $K^4_{NSD}$ will
be connected, and its explicit structure will
reveal the $K^4_{NSD}$-atoms.
Let us assume in the following that for a given $n$ the space
$K^n_{NSD}$ is connected.

For $n\equiv 3 \,\mbox{mod}\, 4$ resp. $n\equiv 1 \,\mbox{mod}\, 4$
corresponding to the $K^n_{NSD}$-atom $\ii(n)$ resp. $\iv(n)$ there are
in any case $2$ atoms of $K^n_{or,NSD}=K^n_+\oplus K^n_-$,
either $\ii(n)^R$ and $\ii(n)^L$, or
resp. $\iv(n)^R$ and $\iv(n)^L$.

Similarly, we could pick for arbitrary $n$ any $K^n_{NSD}$-atom $a$
and will find a corresponding pair of $K^n_{or,NSD}$-atoms $a^R$ and $a^L$.

At this place, let us make
the convention to assign the right atom $a^R$ to $K^n_+$ and
the left atom $a^L$ to $K^n_-$.

Now consider all other pairs of dual points $A^R$ and $A^L$ in
$K^n_{or,NSD}$, which constitute the preimage $\pi^{-1}(A)$ of a
non-selfdual point $A\in K^n_{NSD}$.
For any limit $A\to B$ or $C\to A$ in $K^n$ there
exists a corresponding pair of limits $A^{R/L}\to B'$ or $C'\to A^{R/L}$
in $K^n_{or}$, with $B'\in \pi^{-1}(B)$ resp. $C'\in \pi^{-1}(C)$.
Note however that there are no transitions or limits
between conjugate points,
neither $A^R\to A^L$ nor $A^L\to A^R$, because limits cannot reverse
the orientation.

If $B'$ or $C'$ is non-selfdual,
we demand it, as the limit $B'=B^{R/L}$ resp.
the prelimit $C'=C^{R/L}$ of $A^{R/L}$,
to be contained in the same component
of $K^n_{or,NSD}$ as $A^{R/L}$ itself.

Under consideration of the transitivity of
transitions in $K^n_{or}$ and use of the assumed connectedness
of $K^n_{or,NSD}$,
it follows from assignments for the non-selfdual
atoms made above that, {\em all} right algebras have to be in $K^n_+$ and
{\em all} left algebras have to be in $K^n_-$.

If $K^n_{NSD}$ is connected,
this choice is the only one which makes
each of $K^n_+$ and $K^n_-$ connected
and both disconnected to each other. Therefore
it is the canonical assignment in the case of connected
 $K^n_{NSD}$. This will be
the relevant situation in the following sections.
\nl

For $n\geq 4$ let us define a selfdual algebra $A^a_{n,2}$ with
{Abel}ian ideal $I^{(n-1)}$ by
$\NJNF(A^a_{n,2}):=[a\cdot\NJNF(A_2)] \oplus \NJNF(\iv(n-1))$, where $\oplus$
denotes the direct sum of matrices.

Now it is easy to prove
\nl

{\bf Lemma 5.}
{\em
Within $K^n$ for  $n\geq 3$, the subset $K^n_{SD}$ of selfdual elements in
$K^n_{(or)}$ has the following properties:

If there exists a non-selfdual algebra,
which is the case at least for $n$ odd, then $K^n_{SD}$ is not open.

For $n$ odd  $K^n_{SD}$ is neither open nor closed.
}\nl
{Proof:}
Assume that there exists a non-selfdual algebra;
such an algebra is given by $\iv(n)$ for $n$ odd.
Then there is at least one
$K^n_{NSD}$-atom.
Any $K^n_{NSD}$-atom has a selfdual limit.
Hence, there exists a
selfdual  limit from a non-selfdual sequence
$\imp K^n_{(or),NSD}$ not closed
$\imp K^n_{SD}$ not open.\nl
On the other hand, there are also non-selfdual  limits from selfdual
sequences, like $\VIo\to\II$ for $n=3$
and $A^a_{n,2} \lto \iv(n)$  with $a\to 1$ for odd $n>3$
$ \imp K^n_{SD}$ not closed for odd $n\geq 3$.
\hfill$\Box$ \break
Likewise, each of $K^n_\pm$ is neither open nor closed for $n$ odd.
Note that $K^n_{SD}$ open would imply
$K^n_{SD}=K^n$.

In examination of duality of a given algebra, it is useful
to remind the obvious
\nl

{\bf Lemma 6.}
{\em
For an algebra $A\in K^n$, following assertions are equivalent:
\nl
i) $A$ is selfdual.
\nl
ii) The set $S(A)$ of all subalgebras of $A$ is selfdual.
\nl
iii) The set $J(A)$ of all ideals of $A$ is selfdual.
}
\hfill$\Box$ \break
Note that individual elements of $S(A)$ and $J(A)$ taken for themselves can
be non-selfdual while $A$ is selfdual.

Finally we deal with the case of simple {Lie} algebras.
\nl

{\bf Lemma 7.}
{\em
Simple {Lie} algebras are not selfdual.
}\nl
{Proof:}
A simple $n$-dimensional {Lie} algebra $A_n$ can
be characterized by a
{Cartan-Weyl} basis.
Such a basis consisting of generators
$H_i$, $i=1,\ldots,l=\mbox{rank}A_n$,
which span a maximal {Abel}ian subalgebra
(usually called {Cartan} subalgebra) and $n-l$ generators
$E_\alpha$, each satisfying,
for any nonvanishing generator
$H=\alpha^i H_i$ of the {Cartan} subalgebra,
a root equation
$[H, E_\alpha]=\alpha E_\alpha$ $\ (\ast)$
with root
$\alpha=\alpha^i\alpha_i$.
The commutators
$[E_\alpha,E_\beta]=N_{\alpha\beta}E_{\alpha+\beta}$ $\ (\ast\ast)$
for $\alpha+\beta\neq 0$ and
$[E_\alpha,E_{-\alpha}]=H$ $\ (\ast\ast\ast)$
are nonvanishing.
{}From the root equations $(\ast)$ we see that for any nonvanishing
{Cartan} subalgebra element $H$ (given by its coroots $\alpha^i$)
the reflection $H\to-H$ changes the algebra.
Furthermore by $(\ast\ast)$ and $(\ast\ast\ast)$ also none of the
reflections
$E_\alpha\to-E_\alpha$ keeps the
algebra invariant. Since there is no
reflection keeping
the algebra invariant it can not be selfdual.
\hfill$\Box$ \break

For considerations of the topological structure of $K^3$ and $K^4$
in Sec. 6 and 7 respectively, we will define the notion of parametrical
connectedness
of points in $K^n$
like following:
\nl

{\bf Definition 7.}
$X, Y\in K^n$ are called {\em parametrically connected} iff
there exists a continuous curve
$c: [0,1] \to K^n$ with $c(0)=X$ and $c(1)=Y$ such that,
for all $t_1\leq t_2\in [0,1]$ with
$c(t_1)\neq c(t_2)$, there exist some $t_0\in [t_1,t_2]$
such that  $c(t_1)\neq c(t_0)\neq c(t_2)$.
Otherwise $X, Y\in K^n$ are said to be
{\em parametrically disconnected}.
\hfill $\Box$\break
Note that, in the topology $\kappa^n$,
arcwise connectedness does not imply parametrical connectedness
as defined above.

Furthermore, a set $S\subset K^n$ is called parametrically connected,
iff any two points $X,Y\in S$ are parametrically connected in $S$.

$S\subset K^n$ is a {\em parametrically connected component}
iff $S$ is parametrically connected but not a proper subset of
another parametrically connected set.
\nl
\section{\bf Topology of $K^n$ for $n\leq 4$}
\setcounter{equation}{0}
Sec. 6.1 resumes already existing results for $n\leq 3$,
Sec. 6.2 describes in detail the components and transitions
of $K^4$, and Sec. 6.3 gives some overview over the topological
structure of $K^4$.
\subsection{\bf Structure of $K^n$ for $n\leq 3$}
The {Lie} algebras with $n\leq 3$ are well known and listed, e.g.
by {Patera} and {Winternitz} \PaW.
$K^2$ contains only 2 elements, the {Abel}ian $2A_1$ and $A_2$ represented
by the algebra with $[e_2,e_1]=e_1$ as only nonvanishing bracket.
So $A_2$ has the ideal $J_1=A_1$ spanned by $\{e_1\}$, and is characterized
by $C_{<2>}=(1)\neq 0$ in contrast to $2A_1$. Note that
$A_2\equiv {\V}^{(2)}$. Obviously $\dim K^2_*=0$ and the unimodular
subset $U^2_*\subset K^2_*$ is empty.

The elements of $K^3$ correspond to the famous {Bianchi}
(or {Bianchi-Behr}) types.
They have been classified independently
first by S. {Lie} \Lie\ and then by L. {Bianchi}  \Bi.
For their systematic derivation and explanation of their role for
cosmological models see e.g. {Kramer, Stephani} et al. \Kr.
For convenience of the reader we give
for each of the Bianchi types I up to IX an explicit description
by the commutators of its generators $e_1, e_2, e_3$ according
to {Landau-Lifschitz} \Lan:

Types I,II and VIII/IX are given by basic commutators
$$
[e_1,e_2]=n_3 e_3, [e_2,e_3]=n_1 e_1, [e_3,e_1]=n_2 e_2,
$$
with triplets $(n_1,n_2,n_3)$ respectively given by $(0,0,0), (1,0,0)$ and
$(1,1,\mp 1)$.
The 1-parameter families \VIh/\VIIh with $h\geq 0$ are given respectively
by
$$
[e_1,e_2]=e_3+he_2, [e_2,e_3]=0, [e_1,e_3]=\pm e_2+he_3,
$$
and especially it is $\III={\rm VI}_1$. IV resp. V are given by
$$
[e_1,e_2]=be_3+e_2, [e_2,e_3]=0, [e_1,e_3]=e_3
$$
with $b=1$ resp. $b=0$.

The 3-dimensional real {Lie} algebras in the notation of
{Patera} and {Winternitz} \PaW\ can be
characterized by their NJNF, which is
simultaneously the normal form
(see Eqs. (6.1) up to (6.4) below)
of the {Bianchi} types associated to them like in Table 1.
\bigskip \nl
{\normalsize
\begin{tabular}{ccccccccccc}
$3A_1$&$A_1\oplus A_2$&$A_{3,1}$&$A_{3,2}$&$A_{3,3}$
&$A_{3,4}$&$A^a_{3,5}$&$A_{3,6}$&$A^a_{3,7}$&$A_{3,8}$&$A_{3,9}$\\
I&III&II&IV&V&\VIo&\VIh&\VIIo&\VIIh&VIII&IX
\end{tabular}
\begin{center}
\begin{tabular}{ll}
Table 1: &Inequivalent 3-dim. {Lie} algebras as denoted in \PaW\ (upper row)\\
         &and corresponding {Bianchi} types (lower row).
\end{tabular}
\end{center}
\smallskip }
For convenience we explicitly give the nonvanishing commutators
for the indecomposable algebras $A_{3,1}$ up to $A_{3,9}$ from \PaW:
$$
A_{3,1}: [e_2,e_3]=e_1;
$$
$$
A_{3,2}: [e_1,e_3]=e_1, [e_2,e_3]=e_1+e_2;
$$
$$
A_{3,3}: [e_1,e_3]=e_1, [e_2,e_3]=e_2;
$$
$$
A_{3,4}: [e_1,e_3]=e_1, [e_2,e_3]=-e_2;
$$
$$
A^a_{3,5}: [e_1,e_3]=e_1, [e_2,e_3]=ae_2,\ 0<\vert a\vert<1;
$$
$$
A_{3,6}: [e_1,e_3]=-e_2, [e_2,e_3]=e_1;
$$
$$
A^a_{3,7}: [e_1,e_3]=ae_1-e_2, [e_2,e_3]=e_1+ae_2,\ 0<a;
$$
$$
A_{3,8}: [e_1,e_2]=e_1, [e_2,e_3]=e_3, [e_3,e_1]=2e_2;
$$
$$
A_{3,9}: [e_1,e_2]=e_3, [e_3,e_1]=e_2, [e_2,e_3]=e_1.
$$
In $3$ dimensions, all solvable algebras contain the {Abel}ian ideal
$J_2=2A_1$.
Therefore they can be characterized by their NJNF.
$$
\NJNF(\I) =
\left(
\begin{array}{cc}
 0 &   \\
   & 0
\end{array}
\right) ,
$$
$$
\NJNF(\III) =
\left(
\begin{array}{cc}
 1 &   \\
   & 0
\end{array}
\right) ,
$$
$$
\NJNF(\II) =
\left(
\begin{array}{cc}
 0 & 1 \\
   & 0
\end{array}
\right) ,
\qquad \II={\II}^{(3)} ,
$$
$$
\NJNF(\IV) =
\left(
\begin{array}{cc}
 1 & 1 \\
   & 1
\end{array}
\right) ,
\qquad \IV={\IV}^{(3)} ,
$$
\begin{equation}
\NJNF(\V) =
\left(
\begin{array}{cc}
 1 &  \\
   & 1
\end{array}
\right) ,
\qquad \V={\V}^{(3)} .
\end{equation}
\begin{equation}
\NJNF(\VIo) =
\left(
\begin{array}{cc}
 1 &  \\
   & -1
\end{array}
\right) ,
\
\NJNF(\VIh) =
\left(
\begin{array}{cc}
 1 &  \\
   & a
\end{array}
\right) .
\end{equation}
In Eq. (6.2) the range $0< h< \infty, h\neq 1$ of the parameter
according to {Landau-Lifschitz} \Lan, denoted here by  $h$,
is monotonously homeomorphic to the range
$-1< a< 1, a\neq 0$. $h=1$ resp. $a=0$ yields a decomposable algebra,
namely \III.

\begin{equation}
\NJNF(\VIIo) =
\left(
\begin{array}{cc}
 0  & 1 \\
 -1 & 0
\end{array}
\right) ,
\
\NJNF(\VIIh) =
\left(
\begin{array}{cc}
 a  & 1 \\
 -1 & a
\end{array}
\right) .
\end{equation}
In Eq. (6.3) the range $0< h< \infty$ of the parameter
according to {Landau-Lifschitz} \Lan, denoted here by  $h$,
is monotonously homeomorphic to the range
$0< a< \infty$.\nl
Note that for a topological characterization of $K^3$ it is sufficient to
know the relation of the parameters $a$ and $h$ in (6.2) and (6.3)
at points of qualitative change in the NJNF and to ensure
homeomorphisms of the ranges between these critical points.
This is precisely the data we have given above. (Though the explicit
relation of $a$ and $h$ follows from the equivalence transform to normal form,
here we do not need to calculate it.)
Both \VIh and \VIIh are unimodular for $h=0$ and converge to \II for
$0\leq h<\infty$ and to \IV for $h\to \infty$.\nl
The simple algebras $\VIII=su(1,1)$ and $\IX=su(2)$
are described respectively by the 3 matrices
\begin{equation}
C_{<3>} =
\left(
\begin{array}{cc}
 0  & 1 \\
 -1 & 0
\end{array}
\right) ,
\qquad
C_{<1>} = -C_{<2>} =
\left(
\begin{array}{cc}
 0     &  1 \\
 \pm 1 &  0
\end{array}
\right) .
\end{equation}
$\NJNF(C_{<3>})=\NJNF(\VIIo)$ for both $\VIII$ and $\IX$, \nl
but
$\NJNF(C_{<1>})=\NJNF(C_{<2>})$ is equal to $\NJNF(\VIIo)$ for $\VIII$ and to
$\NJNF(\VIo)$ for $\IX$. Therefore $\IX\to\VIIo$, but $\IX\not\to\VIo$,
but both $\VIII\to\VIo$ and $\VIII\to\VIIo$.  \nl
Considering all components, their parametrical limits and transitions
together,
we get the full topological structure of $K^3$,
which includes a transitive network of nearest neighbour
transitions between different components. The network has been depicted
already by {Mac Callum} \Mac\ and its transitivity
was outlined by {Schmidt} \aSch.
We have $\dim K^3=1$, since its largest parametrically connected components
are 1-dimensional. \nl
For the unimodular subvariety $U^3_*\subset K^3_*$ it is
$\dim U^3_* =0$, and $\{\VIII,\IX\}\subset U^3_*$ is a minimal dense
subset of isolated points.
Fig. 1 shows the transitive network of transitions in $K^3_*$,
with unimodular points encircled.

\vspace*{12.7truecm}
\begin{center}
{\normalsize Fig. 1: Transitive network of transitions in $K^3_*$.}
\end{center}
\np\noindent
\subsection
{\bf Components of $K^4$, transitions and parametrical limits}
The real 4-dimensional {Lie} algebras have been classified
by {Mubarakzyanov} \bMu\ and
listed by {Patera} and {Winternitz} \PaW. An early, somehow more
coarse classification has been given by {Petrov} \Pe.
For the convenience of the reader we explicitly give this classification
in terms of nonvanishing basic commutators. In order to avoid confusion with
the $3$-dimensional {Bianchi} types we alter the notation of {Petrov}'s
classes \Pe\ from I,\ldots,VIII to $\wp_i,\ i=1,\ldots,8$.
The subclasses ${\rm VI}_{1/4}$ are written as \PVIac respectively,
and ${\rm VI}_2$ together with ${\rm VI}_3$ are resumed in a single class
\PVIb in order to correspond to distinct classes of \PaW. With this
notation {Petrov}'s classes are characterized like following:

Solvable algebras, without {Abel}ian subgroup $3A_1$:
$$
\PI: [e_2,e_3]=e_1, [e_1,e_4]=ce_1, [e_2,e_4]=e_2, [e_3,e_4]=(c-1)e_3,\
c\in\R;
$$
$$
\PII: [e_2,e_3]=e_1, [e_1,e_4]=2e_1, [e_2,e_4]=e_2, [e_3,e_4]=e_2+e_3;
$$
$$
\PIII: [e_2,e_3]=e_1, [e_1,e_4]=qe_1, [e_2,e_4]=e_3, [e_3,e_4]=-e_2+qe_3,\
q^2<4;
$$
$$
\PIV: [e_2,e_3]=e_2, [e_1,e_4]=e_1;
$$
$$
\PV: [e_2,e_3]=e_2, [e_3,e_1]=-e_1, [e_1,e_4]=e_2, [e_2,e_4]=-e_1;
$$

Solvable algebras, with {Abel}ian subgroup $3A_1$:
$$
\PVIa: [e_1,e_4]=ae_1+be_4, [e_2,e_4]=ce_2+de_4, [e_3,e_4]=ee_3+fe_4,
$$
with real tuples $(a,b,c,d,e,f)$ of the form
$(0,0,0,0,0,0)$, $(0,1,0,1,0,0)$, $(0,1,0,1,0,1)$, $(1,1,0,0,0,0)$ or
$(1,0,c,0,e,0)$;
$$
\PVIb: [e_1,e_4]=ke_1+e_2, [e_2,e_4]=ke_2+de_3, [e_3,e_4]=ee_3,\
k\in \R,\ d,e\in\{0,1\};
$$
$$
\PVIc: [e_1,e_4]=ke_1+e_2, [e_2,e_4]=-e_1+ke_2, [e_3,e_4]=le_3,\
k,l\in \R;
$$
Non solvable algebras:
$$
\PVII: [e_1,e_2]=e_1, [e_2,e_3]=e_3, [e_3,e_1]=-2e_2;
$$
$$
\PVIII: [e_1,e_2]=e_3, [e_2,e_3]=e_1, [e_3,e_1]=e_2.
$$
\np
In the following we use the characterization of equivalence
classes by their NJNF, according to Sec. 3, in order to find relative
positions of the equivalence classes in $K^4$, possible
transitions between them and parametrical limits of parametrically
connected components of $K^4$.
%
%
\nl\nl
{\bf 6.2.1 Decomposable {Lie} algebras}
\nl\nl
A decomposable 4-dimensional {Lie} algebra can have the structures
$4A_1$, $2A_1\oplus A_2$, $2A_2$ or $A_1\oplus A_3$. The first 3
possibilities are unique, since $A_1$ is the unique 1-dim. {Lie} algebra
and $A_2$ is the unique non {Abel}ian 2-dim. {Lie} algebra. Note that
$2A_2\equiv \PIV$ in {Petrov}'s classification \Pe.
$A_1\oplus A_3$ consists of 9 classes, given by $\{A_{3,i}\}_{i=1,\ldots,9}$
listed
in Table 1. It is $A_1\oplus \II\equiv {\II}^{(4)}$.
$A_1\oplus \VIII$
and $A_1\oplus\IX$ are the same as in \Pe\ the \PVII and \PVIII respectively.
\nl
Transitions and limits: Besides the transitions and limits induced by
{Abel}ian embedding
$\oplus \R$ of transitions in $K^3$, there are further transitions,
which prevent the embedding $\oplus \R$ to be a homeomorphism.
So for example $\V\oplus \R\to \II\oplus \R$, but
$\V\not\to \II$.
This demonstrates that, while $\V$ is an atom, $\V\oplus \R$ is not.
Furthermore $\VIo\oplus \R$ and $\VIIo\oplus \R$ both go first to $A_{4,1}$
and then to $\II\oplus \R$.

$\VIII\oplus \R$ has a limit in the non decomposable
$A_{4,8}$ and, like $\IX\oplus \R$, also in $A_{4,10}$,
as described below.

The algebra $2A_2$ in spite of being decomposable  is not the limit of
any other algebra in $K^4$. It has transitions to $\VIh\oplus \R$
with $h\geq 0$, to $A_{4,3}$ and to $A^0_{4,9}$.
\nl\nl
\noindent
{\bf 6.2.2 Indecomposable {Lie} algebras}
\nl\nl
Coarsely these algebras have already been classified by {Petrov} \Pe.
Table 2
relates his classification to that of {Patera} and
{Winternitz} \PaW.
\smallskip
{\normalsize
\begin{center}
\begin{tabular}{ccccccc}
$A_{4,1..4}$&$A_{4,5}$&$A_{4,6}$&$A_{4,7}$&$A_{4,8/9}$&$A_{4,10/11}$&$A_{4,12}$
\smallskip\\
\PVIb         &\PVIa  &\PVIc  &\PII   &\PI        &\PIII        &\PV
\end{tabular}\smallskip\\
\begin{tabular}{ll}
Table 2:&Classification of {Petrov} \Pe\ (lower row)
and \PaW\ (upper row)\\
        &  of 4-dim. {Lie} algebras except decomposable ones.
\end{tabular}
\end{center}
\smallskip }
The algebra \PIII, $q=0$ is the
same as $A_{4,10}$. It is the only indecomposable $4$-dimensional
algebra that corresponds to a maximal isometry group of a $3$-dimensional
homogeneous {Riemann}ian space
(see Sec. 8, 9 below and {Bona} and {Coll} \Bo, Theorem 1).

For convenience we explicitly give the nonvanishing commutators
of the indecomposable algebras $A_{4,1}$ up to $A_{4,12}$ according
to \PaW:

$$
A_{4,1}: [e_2,e_4]=e_1, [e_3,e_4]=e_2;
$$
$$
A^a_{4,2}: [e_1,e_4]=ae_1, [e_2,e_4]=e_2, [e_3,e_4]=e_2+e_3,\ a\neq 0;
$$
$$
A_{4,3}: [e_1,e_4]=e_1, [e_3,e_4]=e_2;
$$
$$
A_{4,4}: [e_1,e_4]=e_1, [e_2,e_4]=e_1+e_2, [e_3,e_4]=e_2+e_3;
$$
$$
A^{a,b}_{4,5}: [e_1,e_4]=e_1, [e_2,e_4]=ae_2, [e_3,e_4]=be_3,\
-1\leq a\leq b\leq 1,\ ab\neq 0;
$$
$$
A^{a,b}_{4,6}: [e_1,e_4]=ae_1, [e_2,e_4]=be_2-e_3, [e_3,e_4]=e_2+be_3,\
b\geq 0,\ a\neq 0;
$$
$$
A_{4,7}: [e_1,e_4]=2e_1, [e_2,e_4]=e_2, [e_3,e_4]=e_2+e_3, [e_2,e_3]=e_1;
$$
$$
A_{4,8}: [e_2,e_3]=e_1, [e_2,e_4]=e_2, [e_3,e_4]=-e_3;
$$
$$
A^b_{4,9}: [e_2,e_3]=e_1, [e_1,e_4]=(1+b)e_1, [e_2,e_4]=e_2, [e_3,e_4]=be_3,\
-1<b\leq 1;
$$
$A_{4,8}$ is the parametrical limit of $A^b_{4,9}$ for $b\to -1$;
hence by {Mubarakzyanov} \aMu\ and {Petrov} \Pe\
$A_{4,8}$ and $A_{4,9}$ are subsumed in a single
$1$-parameter set.
$$
A_{4,10}: [e_2,e_3]=e_1, [e_2,e_4]=-e_3, [e_3,e_4]=e_2;
$$
$$
A^a_{4,11}: [e_2,e_3]=e_1, [e_1,e_4]=2ae_1, [e_2,e_4]=ae_2-e_3,
[e_3,e_4]=e_2+ae_3,\ 0<a;
$$
$A_{4,10}$ is the parametrical limit of $A^a_{4,11}$ for $a\to 0$;
hence by {Mubarakzyanov} \aMu\ and {Petrov} \Pe\
$A_{4,10}$ and $A_{4,11}$ are subsumed in a single
$1$-parameter set.
$$
A_{4,12}: [e_1,e_3]=e_1, [e_2,e_3]=e_2, [e_1,e_4]=-e_2, [e_2,e_4]=e_1.
$$
The only difference of this classification to that of
{Mubarakzyanov} \aMu\ is that, unlike there,  here
the endpoints $A_{4,8}$ and $A_{4,10}$ are distinguished
against the rest of the $1$-parameter sets $A_{4,9}$ and $A_{4,11}$
respectively.

In the following we reclassify these algebras by their NJNF.
\nl
\nl
{\bf a) Algebras with an {Abel}ian ideal $J_3=3A_1\equiv \I$}
\nl
\nl
These are the algebras of type \PVI. In the following the cases i) and ii)
correspond to \PVIb, case iii) to \PVIa and case iv) to
\PVIc.
\medskip\hfill\break
i) 1 eigenvalue with 1 Jordan block:\nl
Either the eigenvalue is $\lam = 0$ or otherwise it can be normalized to
$\lam = 1$.
\np
\begin{equation}
\NJNF(A_{4,1}) =
\left(
\begin{array}{ccc}
0 & 1 &   \\
  & 0 & 1 \\
  &   & 0
\end{array}
\right) ,
\
\NJNF(A_{4,4}) =
\left(
\begin{array}{ccc}
1 & 1 &   \\
  & 1 & 1 \\
  &   & 1
\end{array}
\right) .
\end{equation}
Transitions: Obviously $A_{4,4}\to A_{4,1}$ and, by increasing the geometric
multiplicity, $A_{4,4}\to{\IV}^{(4)}\equiv A^1_{4,2}$ resp.
$A_{4,1}\to{\II}^{(4)}\equiv \II\oplus \R$.\medskip\hfill\break
ii) Maximally 2 eigenvalues with together 2 Jordan blocks:\nl
Here JNF$(A)$ consists of both a $1\times 1$ and a $2\times 2$ Jordan block,
with eigenvalues $\lam_1$ and $\lam_2$ respectively. If $\lam_1 = 0$, the
algebra would become decomposable ($\II \oplus \R$ if $\lam_2 = 0$,
$\IV \oplus \R$ if $\lam_2 \neq 0$). Therefore assume $\lam_1 = a \neq 0$.
Either $\lam_2 = 0$, then $\lam_1 = 1$  after normalization, or
$\lam_2 \neq 0$, then it can be normalized to $\lam_2 = 1$. If
$\lam_2 = \lam_1 = a$, there is only 1 eigenvalue, which can be
normalized to $a=1$. Note that $A^1_{4,2}\equiv {\IV}^{(4)}$, which is a
case to be considered separately.
\begin{equation}
\NJNF(A_{4,3}) =
\left(
\begin{array}{ccc}
1 &   &   \\
  & 0 & 1 \\
  &   & 0
\end{array}
\right) ,
\
\NJNF(A^a_{4,2}) =
\left(
\begin{array}{ccc}
a &   &   \\
  & 1 & 1 \\
  &   & 1
\end{array}
\right) .
\end{equation}
Transitions and limits: From $A^a_{4,2}$ with $a\neq 0,1$
to $\IV\oplus \R$ for $a\to 0$,
to $A_{4,4}$ for $a\to 1$, and to $A_{4,3}$ for $\vert a\vert\to\infty$. By
increasing the geometric multiplicity, to $A^{1,a}_{4,5}$ for $0<\vert a\vert
<1$,
to $A^{1,-1}_{4,5}=A^{-1,-1}_{4,5}$ for $a=-1$ and to
$A^{\frac{1}{a},\frac{1}{a}}_{4,5}$ for $1<\vert a\vert <\infty$.
Also generally $A^a_{4,2}\to A_{4,1}$. \nl
{}From ${\IV}^{(4)}\equiv A^{1}_{4,2}$ to ${\II}^{(4)}\equiv \II \oplus \R$
and ${\V}^{(4)}\equiv A^{1,1}_{4,5}$, according to the remark at the theorem
in Sec. 4.
\nl
{}From $A_{4,3}$ to $A_{4,1}$ and, by increasing of geometric multiplicity,
to $\III\oplus \R$.\medskip\hfill\break
\noindent
iii) 3 real eigenvalues as Jordan blocks:\nl
Assuming the largest eigenvalue normalized to $\lam_1 = 1$, there remain
$\lam_2=a$ and $\lam_3=b$ with $-1\leq b\leq a\leq 1$. If $a\cdot b = 0$,
the algebra becomes decomposable ($a=b=0$ yields $\III \oplus \R$, for
$a=1,b=0$ it is $\V \oplus \R$, for $a=0,b=-1$ it is $\VIo \oplus \R$ and
otherwise $a=0$ or $b=0$ yields $\VIh \oplus \R$). Therefore assume
$a\cdot b\neq 0$. The case $a=b=1$ (single 3-fold degenerate eigenvalue)
corresponds to the pure vector type
 $A^{1,1}_{4,5} \equiv \Vv \neq \V \oplus \R$.
In $A^{1,b}_{4,5}$ and $A^{a,a}_{4,5}$ there are 2 eigenvalues, one of them
2-fold degenerate.
\np\noindent
For the nondegenerate case, $-1\le b< a< 1$. Note that
$A^{a,b}_{4,5}=A^{b,a}_{4,5}$, since permutations are in ${\GL}(4)$.
\begin{equation}
\NJNF(A^{a,b}_{4,5}) =
\left(
\begin{array}{ccc}
1 &   &   \\
  & a &   \\
  &   & b
\end{array}
\right) .
\end{equation}
Transitions and limits: From $A^{a,b}_{4,5}, a>b,$
to $A^{\frac{1}{a}}_{4,2}$ for ${b\to a}$,
to $A^{b}_{4,2}$ for ${a\to 1}$.
To $A_{4,4}$ for ${{a\to 1}\atop{b\to 1}}$,
to $\IV\oplus \R$ for ${{a\to 1}\atop{b\to 0}}$,
to $\VIIh\oplus \R, 1<h<\infty,$ for ${{a\not\to 0,1}\atop{b\to 0}}$,
to $A_{4,3}$ for ${{a\to 0}\atop{b\to 0}}$,
to $\VIh\oplus \R, 0\leq h<1,$ for ${{a\to 0}\atop{b\not\to 0}}$,
and to $A^{-1}_{4,2}$ for ${b\to -1}$ and ${a\to \pm 1}$.
Also generally $A^{a,b}_{4,5}\to A_{4,1}$. \nl
Note furthermore that $A^{a,-1}_{4,5}=A^{-a,-1}_{4,5}$. \nl
{}From $A^{1,b}_{4,5}$
to ${\IV}^{(4)}\equiv A^{1}_{4,2}$ for ${b\to 1}$, and
to ${\V}\oplus \R$ for ${b\to 0}$.
Generally $A^{1,b}_{4,5}\to {\II}^{(4)}$.\nl
{}From $A^{a,a}_{4,5}$
to ${\IV}^{(4)}\equiv A^{1}_{4,2}$ for ${a\to 1}$, and
to ${\III}\oplus \R$ for ${a\to 0}$.
Generally $A^{a,a}_{4,5}\to {\II}^{(4)}$.\nl
{}From ${\V}^{(4)}\equiv A^{1,1}_{4,5}$ only to $4A_1$, according to the
theorem of Sec. 4.\medskip\hfill\break
iv) 1 real eigenvalue and 2 complex conjugates:\nl
If $\lam_{2,3}=r(\cos\theta \pm i\sin\theta)$, by normalization
$r\sin\theta = 1$ can be achieved, if $\lam_2 \neq \lam_3$ is
assured (otherwise the Jordan block becomes diagonal). Set then
$r\cos\theta = b$ and $\lam_1=a$. Demand $a\neq 0$ to exclude
decomposability ($a = 0$ yields $\VIIh \oplus \R$ and $b=0$ then
corresponds to $h=0$) and without restriction $b\geq 0$.
\begin{equation}
\NJNF(A^{a,b}_{4,6}) =
\left(
\begin{array}{ccc}
a &   &   \\
  & b & 1 \\
  &-1 & b
\end{array}
\right) .
\end{equation}
Transitions and limits: For $a\to 0$, $A^{a,b}_{4,6}\to \VIIh\oplus \R$,
with $0\leq h<\infty$
corresponding to $0\leq b<\infty$.
For a fixed ratio $\frac{a}{b}$ and
$b\to\infty$ there is a limit to $A^{\frac{a}{b}}_{4,2}$, if $a\neq b$,
and to $A_{4,4}$, if $a=b$.
$A^{a,b}_{4,6}\to A_{4,3}$ for $b$ finite (esp. $b=0$) and
$\vert a\vert\to\infty$, and
$A^{a,b}_{4,6}\to \IV\oplus \R$ for $a$ finite (esp. $a=0$) and $b\to\infty$.
Also generally $A^{a,b}_{4,6}\to A_{4,1}$.\nl
Note furthermore that $A^{a,0}_{4,6}=A^{-a,0}_{4,6}$.
\nl
\nl
\np\noindent
{\bf b) Algebras with a nilpotent ideal $J_3=A_{3,1}\equiv{\rm II}$}
\nl
\nl
In the following case i) corresponds to \PII,  case ii) to \PI
and case iii) to \PIII.\medskip\hfill\break
i) 2 eigenvalues with together 2 Jordan blocks:\nl
\begin{equation}
\NJNF(A_{4,7}) =
\left(
\begin{array}{ccc}
2 &   &   \\
  & 1 & 1 \\
  &   & 1
\end{array}
\right) .
\end{equation}
Transitions: $A_{4,7}\to A^{2}_{4,2}$ for $J_3\to \I$.
Furthermore
$A_{4,7}\to A^{1}_{4,9}$.\medskip\hfill\break
ii) 3 real eigenvalues as Jordan blocks:\nl
$$
\NJNF(A^{b}_{4,9}) =
\left(
\begin{array}{ccc}
1+b &   &   \\
    & 1 &   \\
    &   & b
\end{array}
\right) , 0 < \vert b\vert < 1 ,
$$
$$
\NJNF(A^{0}_{4,9}) =
\left(
\begin{array}{ccc}
1 &   &   \\
  & 1 &   \\
  &   & 0
\end{array}
\right) ,
\qquad
\NJNF(A^{1}_{4,9}) =
\left(
\begin{array}{ccc}
2 &   &   \\
  & 1 &   \\
  &   & 1
\end{array}
\right) ,
$$
\begin{equation}
\NJNF(A_{4,8}) =
\left(
\begin{array}{ccc}
0 &   &   \\
  & 1 &   \\
  &   &-1
\end{array}
\right) .
\end{equation}
Transitions: From $A^b_{4,9}$
to $A_{4,8}$ for $b\to -1$,
to $A^0_{4,9}$ for $b\to 0$, and
to $A_{4,7}$ for $b\to 1$.
Furthermore, for $J_3\to \I$, to
$A^{\frac{1}{1+b},\frac{b}{1+b}}_{4,5}$ if $0<b<1$, and to
$A^{{1+b},{b}}_{4,5}$ if $-1<b<0$.\nl

For $J_3\to \I$, $A^1_{4,9}\to A^{\frac{1}{2},\frac{1}{2}}_{4,5}$
and $A_{4,8} \to \VIo\oplus \R$. $A^0_{4,9}$  goes to $\IV\oplus \R$
and further to $\V\oplus \R$.
Since $\VIII\oplus \R\to A_{4,8}$, the latter is a limit from
a decomposable algebra.\medskip\hfill\break
iii) 1 real eigenvalue and 2 complex conjugates:\nl
$$
\NJNF(A^{a}_{4,11}) =
\left(
\begin{array}{ccc}
2a &   &   \\
   & a & 1 \\
   &-1 & a
\end{array}
\right) , a > 0 ,
$$
\begin{equation}
\NJNF(A_{4,10}) =
\left(
\begin{array}{ccc}
0 &   &   \\
  & 0 & 1 \\
  &-1 & 0
\end{array}
\right) .
\end{equation}
\np
Transitions: From $A^a_{4,11}$
to $A_{4,10}$ for $a\to 0$,
to $A^1_{4,9}$ for $a\to \infty$ and,
for $J_3\to \I$, to
$A^{{2a},{a}}_{4,6}$.\nl
For $J_3\to \I$, $A_{4,10}\to \VIIo\oplus \R$. Furthermore
both $\VIII\oplus \R\to A_{4,10}$ and $\IX\oplus \R\to A_{4,10}$.
\nl
\nl
{\bf c) Algebras with a pure vector type ideal $J_3=A_{3,3}\equiv{\rm V}$}
\nl
\nl
This case corresponds to type \PV.
\begin{equation}
\NJNF(A_{4,12}) =
\left(
\begin{array}{ccc}
0 & 1 &   \\
-1& 0 &   \\
  &   & 0
\end{array}
\right) .
\end{equation}
Transitions: $A_{4,12}$ goes to $\V\oplus \R$, to  $\VIIh\oplus \R$,
especially for $J_3\to \I$ to  $\VIIo\oplus \R$, and to $A_{4,9}$.
\nl
\nl

\subsection{\bf The topological structure of $K^4$}
Since we know all components, their parametrical limits and transitions
in $K^4$,
we can now put
them together, in order to determine the full topological structure
of $K^4$.
Fig. 2 a), b) and c)
show components of $K^4_*$, with $J_3$ equal to $\I$, $\II$ and $\V$
respectively, as parts of the transitive network of convergence. The dashed
lines in Fig. 2 a) indicate the $\kappa^4$ limit lines from lines
in Fig. 2 b).
$\dim K^4=2$, since its largest (parametrically connected) components are
2-dimensional.
\nl

For the unimodular subvariety $U^4_*\subset K^4_*$ it is
$\dim U^4_* =1$. The union of
$\VIII\oplus \R$, $\IX\oplus \R$,
$\{A^{a,-a-1}_{4,5},-\frac{1}{2}<a<0\}$
and
$\{A^{-2b,b}_{4,6},0<b<\infty\}$ is a dense
subset of $ U^4_*$ and consists of a minimum number of parametrically
connected components, namely 2 isolated points and 2 isolated line segments.
In Fig. 2
the unimodular lines are dotted, and the unimodular points encircled.

\np
\vspace*{17.5truecm}
\begin{center}
{\normalsize Fig. 2 a: Transitions and limits at components of $K^4_*$ with
ideal \I.}
\end{center}
\vspace*{8.5truecm}
\begin{center}
{\normalsize Fig. 2 b: Transitions at components of $K^4_*$ with ideal \II.}
\end{center}
\vspace{6.5truecm}
\begin{center}
{\normalsize Fig. 2 c: Transitions at components of $K^4_*$ with ideal \V.}
\end{center}

\section{\bf Orientation duality in $K^n_{(or)}$ for $n\leq 4$}
\setcounter{equation}{0}
%
%
%
In this section we
examine in detail all points in $K^n_{or}$ for $n\leq 4$ under the
aspect of orientation duality.
%
%
%

In Sec. 7.1 the topological structure of $K^n_{or}$ for
$n\leq 3$ is analysed by use of the (O)NJNF, thus reproducing the
results listed by {Schmidt} \bSch.
The connected components $K^3_\pm$ are determined explicitly.

Using the same method, Sec. 7.2 analyses the orientation duality structure
of $K^4_{or}$ in detail. Especially we determine the connected components
$K^4_\pm$.

\subsection{\bf Structure of $K^n_{or}$ for $n\leq 3$}
The {Lie} algebras in $K^n$ for $n\leq 3$ have been classified in
Sec. 6.1
using their $n-1$-dimensional ideals and the NJNF.
Their orientation duality has already been listed by
{Schmidt} \bSch.

$K^2$ contains only 2 elements, the {Abel}ian $2A_1$ and $A_2$ represented
by the algebra with $[e_2,e_1]=e_1$ as only nonvanishing bracket. Both are
selfdual, because e.g. $e_1\to-e_1$ does not change the algebra.
So $K^2_{or}=K^2_{SD}=K^2$

The elements of $K^3$ correspond to the familiar Bianchi types.
In the following
we analyse the orientation duality by looking at the NJNF in $K^3$ and for
non-selfduality also considering the ONJNF, defining the elements of
$K^3_\pm$.

The solvable algebras in $K^3$ contain all the {Abel}ian ideal $J_2=2A_1$.
In Sec. 6.1 they are classified according to their NJNF.
Similarly the solvable algebras in $K^3_{or}$ can be classified
according to their ONJNF, which agrees the NJNF in the case of
selfduality. So the selfdual algebras in $K^3_{or}$ correspond to the
following cases of NJNF w.r.t. the {Abel}ian ideal $J_2$:
$$
\NJNF(\I) =
\left(
\begin{array}{cc}
 0 &   \\
   & 0
\end{array}
\right) ,
\quad
\NJNF(\V) =
\left(
\begin{array}{cc}
 1 &  \\
   & 1
\end{array}
\right) ,
$$
$$
\NJNF(\III) =
\left(
\begin{array}{cc}
 1 &   \\
   & 0
\end{array}
\right) ,
$$
\begin{equation}
\NJNF(\VIo) =
\left(
\begin{array}{cc}
 1 &  \\
   & -1
\end{array}
\right) ,
\
\NJNF(\VIh) =
\left(
\begin{array}{cc}
 1 &  \\
   & a
\end{array}
\right) .
\end{equation}
The algebras $\I$ and $\III$ are selfdual, since they are
decomposable. All algebras in Eq. (7.1) invariant under the reflection
$e_1\to-e_1$, which guarantees their selfduality.
The parameter range $0< h< \infty, h\neq 1$
($h$ denoting the parameter of {Landau-Lifschitz} \Lan),
corresponds monotonously to
$-1< a< 1, a\neq 0$.
$h=1$ resp. $a=0$ yields the decomposable \III.
So $K^3_{SD}=\{\I,\V\}\cup\{\VIh,0\leq h<\infty\}$.

The other solvable algebras which are not invariant under any reflection
are non-selfdual.
According to Sec. 3 and Sec. 4 we choose the reflection $e_3\to-e_3$
to characterize them as algebras in $K^3_\pm$, with their ONJNF
respectively given like following:
$$
\ONJNF\{\IIRL\}=\pm\NJNF(\II) =
\pm
\left(
\begin{array}{cc}
 0 & 1 \\
   & 0
\end{array}
\right) ,
$$
$$
\ONJNF\{\IVRL\}=\pm\NJNF(\IV) =
\pm
\left(
\begin{array}{cc}
 1 & 1 \\
   & 1
\end{array}
\right) ,
$$
$$
\ONJNF\{\VIIoRL\}=\pm\NJNF(\VIIo) =
\pm\left(
\begin{array}{cc}
 0  & 1 \\
 -1 & 0
\end{array}
\right) ,
$$
\begin{equation}
\ONJNF\{\VIIhRL\}=\pm\NJNF(\VIIh) =
\pm\left(
\begin{array}{cc}
 a  & 1 \\
 -1 & a
\end{array}
\right) .
\end{equation}
In Eq. (7.2) the parameter range $0< h< \infty$
($h$ denoting the parameter of {Landau-Lifschitz} \Lan)
corresponds monotonously to the range
$0< a< \infty$.

The simple algebras $\VIII=su(1,1)$ and $\IX=su(2)$
are described respectively by the 3 matrices
$$
C_{<3>} =
\left(
\begin{array}{cc}
 0  & 1 \\
 -1 & 0
\end{array}
\right) ,
\qquad
C_{<1>} = -C_{<2>} =
\left(
\begin{array}{cc}
 0     &  1 \\
 \pm 1 &  0
\end{array}
\right) .
$$
$\NJNF(C_{<3>})=\NJNF(\VIIo)$ for both $\VIII$ and $\IX$, \nl
but
$\NJNF(C_{<1>})=\NJNF(C_{<2>})$ is equal to $\NJNF(\VIo)$ for $\VIII$
and to
$\NJNF(\VIIo)$ for $\IX$.

In the {Cartan-Weyl} basis $H:=-ie_3$, $E_{\pm}:=e_1\pm i e_2$
the nonvanishing commutators are given
as $[H,E_\pm]=\pm E_\pm$ and, for \VIII or \IX respectively,
$[E_+,E_-]=\pm 2H$.
Note that the latter are different real sections in the same complex algebra.

According to Lemma 4.5 neither $\VIII$ nor $\IX$ are selfdual.
We discriminate the right and left algebra by the reflection
$e_3\to -e_3$,
defining both  pairs $\VIIIRL$ and $\IXRL$ of points in
$K^3_\pm$. So it is
\begin{equation}
\ONJNF\{C^{R/L}_{<3>}\}=\ONJNF\{\VIIoRL\}
\end{equation}
and $C_{<1>}$ and $C_{<2>}$ interchange under this reflection.
\nl

The table below summarizes the duality properties of the Bianchi classes in
$K^3$. Note that for any point $A \in K^3\sm K_{SD}$ there exists a
pair $(A^R,A^L) \in K_+\oplus K_-$ of points in $K^3_{or}\sm K_{SD}$
with right/left handed bases respectively.
\medskip \nl
{\normalsize
\begin{tabular}{ccccccccccc}
$3A_1$&$A_1\oplus A_2$&$A_{3,1}$&$A_{3,2}$&$A_{3,3}$
&$A_{3,4}$&$A^a_{3,5}$&$A_{3,6}$&$A^a_{3,7}$&$A_{3,8}$&$A_{3,9}$\\
I&III&II&IV&V&\VIo&\VIh&\VIIo&\VIIh&VIII&IX\\
1&1  &0 &0 &1&1   &1   &0    &0    &0   &0
\end{tabular}
\begin{center}
\begin{tabular}{ll}
Table 1: &3-dimensional {Lie} algebra classes in $K^3$,\\
         &corresponding Bianchi types and selfduality (yes=1/no=0)
\end{tabular}
\end{center}
\medskip }

The non-selfdual subset of $K^3_{or}$ has 2 connected $1$-dimensional
components, $K^3_+$ and $K^3_-$, given respectively by
\begin{equation}
\VIIIRL/\IXRL\to\VIIoRL\gets\VIIhRL\to\IVRL\to\IIRL ,
\end{equation}
where $\IIRL$ is respectively the atom of $K^3_\pm$.
\subsection{\bf Structure of $K^4_{or}$}
In Sec. 6.2 we classified the real 4-dimensional {Lie} algebras.
In this section they are reconsidered under the aspect
of orientation duality.
\nl\nl
{\bf 7.2.1 Selfdual {Lie} algebras}
\nl\nl
There exist following types of selfdual algebras:
a) all decomposable ones, b) indecomposable ones with ideal \I,
and c) some indecomposable ones with ideal \II.
\nl
\nl
{\bf a) Decomposable ones:}
\nl
\nl
All decomposable {Lie} algebras are selfdual.
A decomposable 4-dimensional {Lie} algebra can have the structures
$4A_1$, $2A_1\oplus A_2$, $2A_2$ or $A_1\oplus A_3$. The first 3
possibilities are unique, since $A_1$ is the unique 1-dim. {Lie} algebra
and $A_2$ is the unique non{Abel}ian 2-dim. {Lie} algebra.
$A_1\oplus A_3$ consists of 9 classes, given by $\{A_{3,i}\}_{i=1,\ldots,9}$
listed in Table 1. Note that $A_1\oplus \II\equiv {\II}^{(4)}$.
\nl
\nl
\np\noindent
{\bf b) Indecomposable ones with 
ideal $J_3={\rm I}$:}
\nl
\nl
Algebras with
{Abel}ian ideal $J_3=3A_1\equiv \I$ are selfdual.
They
are given by the following cases:
\nl
i) 1 Jordan block:
\nl
These algebras are invariant under a combination
of the $3$ reflections $e_i\to -e_i$, $i=1,\ldots,3$.

Either the eigenvalue is $\lam = 0$ or otherwise it can be normalized to
$\lam = 1$.
$$
\NJNF(A_{4,1}) =
\left(
\begin{array}{ccc}
0 & 1 &   \\
  & 0 & 1 \\
  &   & 0
\end{array}
\right) ,
$$
\begin{equation}
\NJNF(A_{4,4}) =
\left(
\begin{array}{ccc}
1 & 1 &   \\
  & 1 & 1 \\
  &   & 1
\end{array}
\right) .
\end{equation}
These algebras are a $4$-dimensional analogue to \II and \IV.
While the latter are non-selfdual their even dimensional analogues
are selfdual.
These algebras are the essential dimensional ones, introduced
in Sec. 5 and denoted by $\ii(4)$ and $\iv(4)$.
$\ii(4)$ is an essential dimensional atom.
\nl
\nl
ii) 2 Jordan blocks:
\nl
All these algebras are all invariant under the reflection $e_1\to -e_1$.

\begin{equation}
\NJNF(A_{4,3}) =
\left(
\begin{array}{ccc}
1 &   &   \\
  & 0 & 1 \\
  &   & 0
\end{array}
\right) ,
\
\NJNF(A^a_{4,2}) =
\left(
\begin{array}{ccc}
a &   &   \\
  & 1 & 1 \\
  &   & 1
\end{array}
\right) .
\end{equation}
In the latter case $a\neq 0$ and $A^1_{4,2}\equiv \IV^{(4)}$.
\nl
\nl
iii) 3 real eigenvalues as Jordan blocks:
\nl
All these algebras are all invariant under the reflection $e_1\to -e_1$.

Assuming the largest eigenvalue normalized to $\lam_1 = 1$, there remain
$\lam_2=a$ and $\lam_3=b$ with $-1\leq b\leq a\leq 1$. If $a\cdot b = 0$,
the algebra becomes decomposable ($a=b=0$ yields $\III \oplus \R$, for
$a=1,b=0$ it is $\V \oplus \R$, for $a=0,b=-1$ it is $\VIo \oplus \R$ and
otherwise $a=0$ or $b=0$ yields $\VIh \oplus \R$). Therefore assume
$a\cdot b\neq 0$. The case $a=b=1$ (single 3-fold degenerate eigenvalue)
corresponds to the pure vector type
 $A^{1,1}_{4,5} \equiv \Vv \neq \V \oplus \R$.
In $A^{1,b}_{4,5}$ and $A^{a,a}_{4,5}$ there are 2 eigenvalues, one of them
2-fold degenerate. For the nondegenerate case, $-1\le b< a< 1$. Note that
$A^{a,b}_{4,5}=A^{b,a}_{4,5}$, since permutations are in ${\GL}(4)$.
\begin{equation}
\NJNF(A^{a,b}_{4,5}) =
\left(
\begin{array}{ccc}
1 &   &   \\
  & a &   \\
  &   & b
\end{array}
\right) .
\end{equation}
\noindent
iv) 1 real eigenvalue and 2 complex conjugates:
\nl
All these algebras are all invariant under the reflection $e_1\to -e_1$.

If $\lam_{2,3}=r(\cos\theta \pm i\sin\theta)$, by normalization
$r\sin\theta = 1$ can be achieved, if $\lam_2 \neq \lam_3$ is
assured (otherwise the Jordan block becomes diagonal). Set then
$r\cos\theta = b$ and $\lam_1=a$. Demand $a\neq 0$ to exclude
decomposability ($a = 0$ yields $\VIIh \oplus \R$ and $b=0$ then
corresponds to $h=0$) and without restriction $b\geq 0$.
\begin{equation}
\NJNF(A^{a,b}_{4,6}) =
\left(
\begin{array}{ccc}
a &   &   \\
  & b & 1 \\
  &-1 & b
\end{array}
\right) .
\end{equation}
\nl
\nl
{\bf c) Indecomposable ones with ideal $J_3={\rm II}$:}
\nl
\nl
There exist algebras with non-selfdual ideal \II, which
are selfdual.
\begin{equation}
\ONJNF(A_{4,8})=\NJNF(A_{4,8}) =
\left(
\begin{array}{ccc}
0 &   &   \\
  & 1 &   \\
  &   &-1
\end{array}
\right).
\end{equation}
This algebra is left invariant by a combination of reflections
$e_4\to-e_4$, $e_1\to-e_1$ and $e_2\leftrightarrow e_3$.
\begin{equation}
\ONJNF(A_{4,10})=\NJNF(A_{4,10}) =
\pm\left(
\begin{array}{ccc}
0 &   &   \\
  & 0 & 1 \\
  &-1 & 0
\end{array}
\right) .
\end{equation}
This algebra is left invariant by a combination of reflections
$e_4\to-e_4$, $e_1\to-e_1$ and $e_2\to -e_2$.
\np\noindent
{\bf 7.2.2 Non-selfdual {Lie} algebras}
\nl\nl
This kind of algebras exists with a basic ideal $J_3$, given either by the
non-selfdual \II or by the selfdual \V.
For all of them we have dual pairs of right and left points in $K^n_{or}$,
which transform to each other by $e_4\to -e_4$, constituting
by Sec. 3 and Sec. 4 the connected components $K^4_\pm$ respectively.
\nl
\nl
{\bf a) Indecomposable ones with ideal $J_3={\rm II}$:}
\nl
\nl
The ideal \II is non-selfdual.
For an algebra $A$ of the kinds listed below the there
exists no reflection leaving the set $J(A)$
invariant.
$$
\ONJNF(A^{R/L}_{4,7})=\pm\NJNF(A_{4,7}) =
\pm\left(
\begin{array}{ccc}
2 &   &   \\
  & 1 & 1 \\
  &   & 1
\end{array}
\right) ,
$$
$$
\ONJNF(A^{b,R/L}_{4,9})=\pm\NJNF(A^{b}_{4,9}) =
\pm\left(
\begin{array}{ccc}
1+b &   &   \\
    & 1 &   \\
    &   & b
\end{array}
\right) , 0 < \vert b\vert < 1 ,
$$
$$
\ONJNF(A^{0,R/L}_{4,9})=\pm\NJNF(A^{0}_{4,9}) =
\pm\left(
\begin{array}{ccc}
1 &   &   \\
  & 1 &   \\
  &   & 0
\end{array}
\right) ,
$$
$$
\ONJNF(A^{1,R/L}_{4,9})=\pm\NJNF(A^{1}_{4,9}) =
\pm\left(
\begin{array}{ccc}
2 &   &   \\
  & 1 &   \\
  &   & 1
\end{array}
\right) ,
$$
\begin{equation}
\ONJNF(A^{a,R/L}_{4,11})=\pm\NJNF(A^{a}_{4,11}) =
\pm\left(
\begin{array}{ccc}
2a &   &   \\
   & a & 1 \\
   &-1 & a
\end{array}
\right) , a > 0 .
\end{equation}
\nl
\nl
{\bf b) Indecomposable ones with ideal $J_3={\rm V}$:}
\nl
\nl
The only case here is given by
\begin{equation}
\ONJNF(A^{R/L}_{4,12})=\pm\NJNF(A_{4,12}) =
\pm\left(
\begin{array}{ccc}
0 & 1 &   \\
-1& 0 &   \\
  &   & 0
\end{array}
\right) .
\end{equation}
Note that besides the selfdual ideal \V there is a second ideal \VIIo
which is not selfdual, causing here the subset of ideals
$S(A_{4,12})$ to be non-selfdual.
Hence $A_{4,12}$ itself is non-selfdual.
\nl
\nl
{\bf c) The space $K^4_{NSD}$ and its components $K^4_\pm$:}
\nl
\nl
Collecting the algebras of the previous subsections a) and b) and recalling
transitions and parametrical limits of components
in $K^4_{NSD}$ according to Sec. 6, we find that
$K^4_{NSD}$ is connected, and so is each of  $K^4_\pm$.
There are $2$ pairs of $K^4_{or,NSD}$-atoms, given by
$A^{0,R/L}_{4,9}$ and $A^{1,R/L}_{4,9}$.
In $n=4$ all $K^4_{NSD}$-atoms have an ideal $\II$, and hence
in the complement of the subspace $K^4_{NSD\vert\I}$ of
$K^4_{NSD}$-algebras with ideal $\I$.

Let us assign e.g. $A^{0,R/L}_{4,9}$ to $K^4_{\pm}$
respectively. Then the connectedness of  $K^4_{\pm}$ and
the orientation preservation of limits within $K^4_{or,NSD}$
imply the assignment $A^{R/L}$ to  $K^4_{\pm}$ respectively.
Note that with these assignments
the component $K^4_+$ is given as
\begin{equation}
\begin{array}{lcrr}
&A^{R}_{4,12}&  \\
&\downarrow&  \\
A^{-1<b<0,R}_{4,9}\to\!&A^{0,R}_{4,9}&\!
\gets A^{0<b<1,R}_{4,9}\to A^{R}_{4,7}\to A^{1,R}_{4,9}\gets A^{a>0,R}_{4,11}
\end{array}
\end{equation}
and the component $K^4_-$ as
\begin{equation}
\begin{array}{lcrr}
&A^{L}_{4,12}&  \\
&\downarrow&  \\
A^{-1<b<0,L}_{4,9}\to\!&A^{0,L}_{4,9}&\!
\gets A^{0<b<1,L}_{4,9}\to A^{L}_{4,7}\to A^{1,L}_{4,9}
\gets A^{a>0,L}_{4,11}
\end{array}
\end{equation}
So the non-selfdual components of $K^4_{(or)}$ are $1$-dimensional.
Note that $A^{R/L}_{4,1}\equiv \ii(4)^{R/L}$ is the atom of $K^4_\pm$
respectively.

\section{\bf Discussion and outlook}
\setcounter{equation}{0}
In Sec. 6 we determined {Lie} algebra transitions
in $K^4$ as limits induced by the topology $\kappa^4$.
Any {In\"on\"u-Wigner} contraction corresponds
to a certain transition; explicitly any of the
{In\"on\"u-Wigner} contractions listed in the tables of
{Huddleston} \Hud\ for real $4$-dimensional {Lie} algebras
corresponds to a transition in $K^4$.
Since {In\"on\"u-Wigner} contractions are only a special
case of the more general {Saletan} contractions,
and since even the latter do not induce all possible transitions
in $K^n$ with $n\geq 3$, it should not be surprising that we
have obtained transitions, which do not correspond to any
{In\"on\"u-Wigner} contraction, like e.g.
transitions $\IX\oplus \R\to A_{4,10}$, $\VIII\oplus \R\to A_{4,10}$
and transitions from $A^{a,b}_{4,5}$, $A^{a,b}_{4,6}$, $\VIh\oplus \R$
and $\VIIh\oplus \R$ to $A_{4,1}$.
The transition $\IX\oplus \R\to A_{4,10}$ corresponds to a
{Lie} algebra contraction, which was given already in \Sal\
(see Eqs. (35') to (37)) as an example of a {Saletan}
contraction, which can not be obtained as
a {In\"on\"u-Wigner} contraction.

It is also interesting to consider transitions in $K^3$ as obtained
in Sec. 5. The {In\"on\"u-Wigner}
contractions for real {Lie} algebras of dimension $d\leq 3$
are classified already by {Conatser} \Co.
The sequence of transitions
$\VIII\to\VIo\to\II\to\I$ is generated
by an iterated {Saletan} contraction (see \Sal, Eqs. (30) and (31)),
applied first to the
{Lie} algebra $\VIII$, of the $3$-dimensional
homogenous {Lorentz} group. On the 4-point subset
$\{\VIII,\VIo,\II,\I\}$ {Saletan} contractions are transitive.
However this transitivity does not hold for {Saletan} transitions
on general subsets of $K^3$.
The sequence of transitions
$\IX\to\VIIo\to\II\to\I$,
starting from the {Lie} algebra $\IX$ of the 3-dimensional
{Euclid}ean group, can not be obtained
by {Saletan} contractions.
The only {Saletan} contractions starting from $\IX$ are in fact
given by a {In\"on\"u-Wigner} contraction
$\IX\to \VIIo$ and the trivial contraction
$\IX\to \I$. Though there exists a different
{In\"on\"u-Wigner} contraction
corresponding to the transitions $\VIIo\to\II$ there is no
{Saletan} contraction corresponding to $\IX\to\II$
(for a proof see \Sal). This example shows that, on an arbitrary
subset of $K^n$ with $n\geq 3$, in general not every transition
can be obtained from a {Saletan} contraction.
It implies that, even
on a set of points connected by {In\"on\"u-Wigner} contractions,
neither {Saletan}
contractions nor {In\"on\"u-Wigner} contractions need to be transitive.

Since we consider transitions between different points in $K^n$,
improper contractions of an algebra to an equivalent one can not be
seen by our method. For $n=4$ {Huddleston} \Hud\ identified
two types of algebras which admit only trivial and improper contractions.
These are precisely the two atoms of $K^4$, namely the unimodular
${\II}^{(4)}\equiv\II\oplus \R$ and
the pure vector type ${\V}^{(4)}\equiv A^{1,1}_{4,5}$.
For arbitrary dimension $n$, the atoms
of $K^n$ have been introduced and described first by
{Schmidt} {\bSch}.

By now the topological properties of $K^n$ for $n\leq 4$ have been examined.
It is natural to demand an investigation for arbitrary dimension $n$.
Practically, this is obstructed by the rapidly increasing number of
equivalence classes for increasing $n$.
A classification
for all nilpotent algebras has been done for $n=6$ by {Morozov} \Mo\
and for $n=7$ by {Ancochea-Bermudez} and {Goze} \An\ in the
complex case and by {Romdhani} \Ro, who
distinguishes $132$ components of real indecomposable nilpotent
$7$-dimensional {Lie} algebras.
For $K^5$ a full classification of all
real {Lie} algebras still distinguishes $40$ components
(compare {Mubarakzjanov} \bMu\ and {Patera} et al. \PaSWZ).
A determination of all possible transitions would be a rather tidy work.
However it is known by  {\PaSWZ}
that $\dim K^5=3$,
because the maximal dimension of its components is 3.
Unlike for the classification of subalgebra
structures of each class in
$K^n$ (see {Patera, Winternitz} \PaW, and {Grigore, Popp} \Gri),
for the determination
of all equivalence classes and transitions between them there exists
no algorithm at present. However, a systematic exploitation of the NJNF,
which has been defined for arbitrary $n$, may contribute
some part to further progress.

The NJNF has proven to be a useful tool in characterizing
distinct $n$-dimensional {Lie} algebras with a common ideal
$J_{n-1}$ as endomorphisms ad$e_n$ of a complementary generator
$e_n$ on that ideal, with characteristic Jordan blocks
of their eigenvalues normalized by an overall scale.
In $4$ dimensions, besides decomposable
algebras, only cases with ideal \I, \II, or \V appear.

In $4$ dimensions, there are no simple algebras.
In general for $n\geq 6$
further classes of simple {Lie} algebras arise, which
lead to an additional further sophistication, as compared to $n=3$.

A combination of the established knowledge on semisimple {Lie} algebras
with the full classification of all {Lie} algebras would be desirable,
but is practically far away, since the dimensionality of the simple
{Lie} algebras increases rapidly with their rank. Note that all
simple components belong to the unimodular subset $U^n_*$.
In Sec. 7 we found that simple {Lie} algebras are non-selfdual
w.r.t. orientation reflection.

In general, we have neither a formula for $\dim K^n$ nor for
$\dim U^n_*$.
%
The $T_0$ topology allows components of different dimensions
to converge pointwise to each other, i.e. such that any point of the first
component converges to some point of the second component and any point of
the second component is the limit of some point of the first.
%

We have determined the topology of the space $K^n_{or}$ for $n\leq 4$.
The essential difference to $K^n$ is that the single non-selfdual
component of the latter is doubled to two  components $K^n_+$
and $K^n_-$.

For $n=3$ or $4$, the space $K^n_{NSD}$ is nonvoid and connected.
For $n=3$ there is a unique $K^3_{NSD}$-atom
$\ii(3)=\II$.
$K^4_{NSD}$ has two atoms, $A^0_{4,9}$ and $A^1_{4,9}$,
and $A^{0<b<1}_{4,9}$ has boundary limits to both of them.

If $K^n_{NSD}$ is connected,
the arbitrariness in assigning  conjugate pairs of points to $K^n_\pm$
can be reduced to a single decision for one pair only,
if we demand that both of $K^n_\pm$ are connected and to each other
disconnected.

At present, for general $n\geq 5$ it is not known whether
$K_{NSD}$ is connected.

{}From Eqs. (7.4) and (7.13-14), we see that the non-selfdual subset
$K^n_{NSD(,or)}$ of $K^n_{(or)}$ is $1$-dimensional
for both, $n=3$ and $n=4$.
We have $\dim K^3_{(or)}=1$ and $\dim K^n_{(or)}\geq 2$ for
dimension  $n\geq 4$: In the latter case the
contribution of the non-selfdual subset is
of dimension less than that of the highest-dimensional
component,
while in the former case it is of highest dimension. Actually,  the  topology
of the highest-dimensional component of $K^3_{or}$ differs
from that of $K^3$ essentially.
The question of dimensionality for general $n\geq 5$ remains open,
for the non-selfdual subset as well as for $K^n_{(or)}$ itself.

Partial progress has been made by determining a candidate of an atom
of the non-selfdual
subset in odd dimension $n$. We found an interesting
periodicity in the structure of this  $K^n_{NSD}$-atom: it is
$\ii(n)$ for $n=3 \,\mbox{mod}\, 4$
and $\iv(n)$ for $n=1 \,\mbox{mod}\, 4$.
However it remains an open problem to determine all atoms of
$K^n_{NSD}$ for arbitrary $n\geq 5$. Presently we do not know
how an atom for even $n\geq 6$ looks like in general.

With Definition 6 the notion of an atom from
{Schmidt} \bSch\ has been generalized to arbitrary subsets.

Although the present work is on the case of real {Lie} algebras,
we want to  make some comments on the analogous complex cases
to the pairs of algebras $\VIIh/\VIh$ ($h\geq 0$),  $\IX/\VIII$ of $K^3$ and
$A^{a,b}_{4,6}/A^{a,b}_{4,5}$, $A^{a}_{4,11}/A^{b}_{4,9}$ and
$A_{4,10}/A_{4,8}$ of $K^4$.
If one considers the analogous $3$- or $4$-dimensional
{Lie} algebras over the
complex basic field the group  $\GL(n)$ is now correspondingly the
group of nonsingular complex linear transformation.
The pairs of complex conjugated eigenvalues associated to the
$2\times 2$ Jordan block of each of the first algebras of the pairs above
in the complex remain as $2$ Jordan blocks in a corresponding
complex JNF. After introducing a similar normalization convention
like for the real case, the complex analogues of the
NJNF will be the same for members of any pair above.
(For $n=3$ this had  already been realized by
{Bianchi} \Bi. The complex $4$-dimensional case was considered
already by {Lie} \Lie).
\nl\nl
{\Large {\bf Acknowledgments}}
\nl\nl
I would like to express my gratitude to H.-J.
{ Schmidt} for
valuable discussions on the present topic.
\np\noindent
{\Large {\bf References}}
\nl\nl
\An\ { J. M. Ancochea-Bermudez} and { M. Goze}.
Classification des algebres de { Lie} nilpotentes de dimension 7.
Arch. Math. {\bf 52}, 175-185 (1989).
\smallskip
\nl
-, Sur la classification des algebres de { Lie} nilpotentes de dimension 7.
Comptes Rendus Acad. Sci. Paris, Ser. I {\bf 302}, 611-613 (1986).
\smallskip
\nl
\Bi\  { L. Bianchi},
Sugli spazii a tre dimensioni che ammettono un gruppo continuo di
movimenti.
Soc. Ital. Sci. Mem. Mat., Ser. {IIIa}, {\bf 11}, 267 (1897).
\smallskip
\nl
\bBi\ -, 
Lezioni sulla teoria dei gruppi continui finiti di trasformazioni,
Cap. XIII: Applicazioni alla teoria degli spazii pluridimensionali
con un gruppo continuo di movimenti,
\S\S 198, 199 (for a description of types \I-\VII and \VIII/\IX)
and \S\S 200-207 (for applications to homogeneous spaces), p. 550-578.
Spoerri, Pisa 1918.
\smallskip
\nl
\Cha\  { R. Charles},
Sur la structure des algebres de { Lie} rigides.
Ann. Inst. Fourier Grenoble {\bf 33}, 3, 65-82 (1984).
\smallskip
\nl
\Ch\  { R. Charles} and { Y. Diakite},
Sur les varietes d'algebres de { Lie} de dimension $\le 7$.
J. Algebra {\bf 91}, 53-63 (1984).
\smallskip
\nl
\Co\  { C. W. Conatser},
Contractions of low-dimensional real { Lie} algebras.
J. Math. Phys. {\bf 13}, 196-203 (1972).
\smallskip
\nl
\Est\  { F. B. Estabrook, H. D. Wahlquist} and { C. G. Behr},
Dyadic analysis of spatially homogeneous world models.
J. Math. Phys. {\bf 9}, 497-504 (1968).
\smallskip
\nl
\Gri\  { D. R. Grigore} and { O. T. Popp},
On the classification of { Lie} subalgebras and { Lie} subgroups.
J. Math. Phys. {\bf 32}, 33-39 (1991).
\smallskip
\nl
\Gru\  { F. Grunewald} and { J. O'Halloran},
Varieties of nilpotent lie algebras of dimension less than six.
J. Algebra {\bf 112}, 315-325 (1988).
\smallskip
\nl
\Hud\ { P. L. Huddleston},
In\"on\"u-Wigner contractions of the real four-dimensional { Lie} algebras.
J. Math. Phys. {\bf 19}, 1645-1649 (1978).
\smallskip
\nl
\In\  { E. In\"on\"u} and { E. P. Wigner},
On the contraction of groups and their representations.
Proc. Nat. Acad. Sci. USA {\bf 39}, 510-524 (1953).
\smallskip
\nl
-, On a particular type of convergence to a singular matrix.
Proc. Nat. Acad. Sci. USA {\bf 40}, 119-121 (1954).
\smallskip
\nl
\Ki\  { A. A. Kirillov} and { Yu. A. Neretin},
The variety $A_n$ of $n$-dimensional { Lie} algebra structures.
Transl., Ser. (2), Amer. Math. Soc. {\bf 137}, 21-30 (1987);
translation from: Some questions in modern analysis, Work. Collect.,
Moskva 1984, 42-56 (1984).
\smallskip
\nl
\np\noindent
\Kr\  { D. Kramer, H. Stephani, M. Mac Callum} and { E. Herlt},
Exact Solutions of Einstein's Field Equations.
Wissenschaften, Berlin 1980.
\smallskip
\nl
\Lan\ { L. D. Landau} and { E. M. Lifschitz},
Lehrb. d. theor. Phys., Bd. II,
12. Aufl., pp. 457 ff.
Akademie-Verlag, Berlin 1992.
\smallskip
\nl
\Lee\  { H. C. Lee},
Sur les groupes de { Lie} r\'eels \`a trois param\`etres.
J. Math. Pures et Appl. {\bf 26}, 251-267 (1947).
\smallskip
\nl
\Lie\  { S. Lie}, Differentialgleichungen, Kap. 21. Chelsea,
Leipzig 1891.
\smallskip
\nl
\bLie\  -, 
Theorie der Transformationsgruppen, Bd. III, Abtheilung VI,
Kap. 28: Allgemeines \"uber die Zusammensetzung $r$-gliedriger
Gruppen, \S\S 136-137, p. 713-732.
Teubner, Leipzig 1893.
\smallskip
\nl
\Mac\  { M. A. H. Mac Callum},
A class of homogeneous cosmological models III: Asymptotic behaviour.
Commun. Math. Phys. {\bf 20}, 57-84 (1971).
\smallskip
\nl
\Mag\  { L. M. Magnin},
Sur les algebres de { Lie} nilpotent de dim $\leq 7$.
J. G. P. {\bf 3}, 1, 119-134 (1986).
\smallskip
\nl
\Mo\  { V. V. Morozov},
Classification of nilpotent { Lie} algebras of 6$^{th}$ order.
Izv. Vyssh. Uchebn. Zaved. Mat. {\bf 4}, 161-171 (1958).
\smallskip
\nl
\aMu\  { G. M. Mubarakzyanov}, On solvable { Lie} algebras. (Russian)
Izv. Vys\v s. U\v cebn. Zavedeni\u\i\ Mat. 1963, no. 1 (32), 114-123.
\smallskip
\nl
\bMu\  -, Classification of real structures of { Lie}
algebras of fifth order. (Russian)
Izv. Vys\v s. U\v cebn. Zavedeni\u\i\ Mat. 1963, no. 3 (34), 99-106.
\smallskip
\nl
\cMu\  -, Classification of solvable { Lie} algebras
of sixth order with a non-nilpotent basis element. (Russian)
Izv. Vys\v s. U\v cebn. Zavedeni\u\i\ Mat. 1963, no. 4 (35), 104-116 (1963).
\smallskip
\nl
\Ne\  { Yu. A. Neretin},
Estimate of the number of parameters assigning an $n$-dimensional algebra.
(Russian)
Isv. Acad. Nauk SSSR, Ser. Mat. {\bf 51}, 306-318 (1987).
\smallskip
\nl
\PaSWZ\ { J. Patera, R. T. Sharp, P. Winternitz} and { H. Zassenhaus},
Invariants of real low dimension { Lie} algebras.
J. Math.  Phys. {\bf 17}, 986-994 (1976).
\smallskip
\nl
\PaW\  { J. Patera} and { P. Winternitz},
Subalgebras of real three- and four-dimensional { Lie} algebras.
J. Math. Phys. {\bf 18}, 1449-1455 (1977).
\smallskip
\nl
\Pe\  { A. S. Petrov}, { Einstein}r\"aume.
Akademie Verlag, Berlin 1964.
\smallskip
\nl
\aRa\  { M. Rainer}, Topology of the space of 4-dimensional real
{ Lie} Algebras. Preprint 93/2, Math. Inst. Univ. Potsdam (1993).
\smallskip
\nl
\bRa\  { M. Rainer}, Orientation Duality of 4-dimensional Real
{ Lie} Algebras.
Preprint 93/8, Math. Inst. Univ. Potsdam (1993).
\smallskip
\nl
\Ri\  { W. Rinow}, Topologie.
Wissenschaften, Berlin 1975.
\smallskip
\nl
\Ro\ { M. Romdhani},
Classification of real and complex nilpotent { Lie} algebras of dimension 7.
Linear Multilin. Algebra {\bf 24}, 167-189 (1989).
\smallskip
\nl
\Sal\ { E. I. Saletan},
Contraction of { Lie} groups.
J. Math. Phys. {\bf 2}, 1-21 (1961).
\smallskip
\nl
\San\ { R. M. Santilli}, Lie-isotopic lifting of Lie symmetries,
I: General considerations. Hadronic J. {\bf 8} 25-35 (1985).
\smallskip
\nl
-, Lie-isotopic lifting of Lie symmetries,
II: Lifting of Rotations. Hadronic J. {\bf 8} 36-51 (1985).
\smallskip
\nl
\aSch\ { H.-J. Schmidt},
Inhomogeneous cosmological models containing homogeneous inner hypersurface
geometry. Changes of the { Bianchi} type.
Astron. Nachr. {\bf 303}, 227-230 (1982).
\smallskip
\nl
\bSch\ -,
On Geroch's limit of space-times and its relation to a new topology in the
space of { Lie} groups.
J. Math. Phys. {\bf 28}, 1928-1936 (1987);
and Addendum in J. Math. Phys. {\bf 29}, 1264 (1988).
\smallskip
\nl
\Se\ { I. E. Segal},
A class of operator algebras which are determined by groups.
Duke Math. J. {\bf 18}, 221-265 (1951).
\smallskip
\nl
\Smrz\ { P. K. Smrz},
Relativity and deformed { Lie} groups.
J. Math. Phys. {\bf 19}, 2085-2088 (1978).
\smallskip
\nl
\aTur\  { P. Turkowski}, Low-dimensional { Lie} algebras.
J. Math. Phys. {\bf 29}, 2139-2144 (1988).
\smallskip
\nl
\bTur\  -, 
Solvable { Lie} algebras of dimension six.
J. Math. Phys. {\bf 31}, 1344-1350 (1990).
\smallskip
\nl
\cTur\  -, 
Structure of real { Lie} algebras.
Linear Algebra Appl. {\bf 171}, 197-212 (1992).
\smallskip
\nl
\Um\  { K. A. Umlauf}, \"uber die Zusammensetzung der endlichen
continuierlichen Transformationsgruppen, insbesondere der Gruppen
vom Range Null.
Breitkopf u. H\"artel, Leipzig 1891.
\smallskip
\nl
\Vra\  { G. Vranceanu}, Lecons de geom\'etrie diff\'erentielle,
p. 105-111. Bucarest 1947.
\smallskip
\nl

\end{document}